\documentclass[leqno,11pt]{article}

\usepackage{amssymb,amsmath,amscd,amsfonts,a4,psfrag,graphicx,latexsym,epsfig,color}
\usepackage[all]{xy}

\newtheorem{theo}{Theorem}[section]
\newtheorem{defi}[theo]{Definition}

\newtheorem{lem}[theo]{Lemma}

\newtheorem{prop}[theo]{Proposition}

\newtheorem{rem}[theo]{Remark}
\newtheorem{coro}[theo]{Corollary}

\newtheorem{exam}[theo]{Example}

\newcommand{\agot}{\ensuremath{\mathfrak{a}}}
\newcommand{\bgot}{\ensuremath{\mathfrak{b}}}

\newcommand{\ugot}{\ensuremath{\mathfrak{u}}}
\newcommand{\kgot}{\ensuremath{\mathfrak{k}}}
\newcommand{\hgot}{\ensuremath{\mathfrak{h}}}
\newcommand{\ggot}{\ensuremath{\mathfrak{g}}}

\newcommand{\tgot}{\ensuremath{\mathfrak{t}}}
\newcommand{\ngot}{\ensuremath{\mathfrak{n}}}
\newcommand{\pgot}{\ensuremath{\mathfrak{p}}}
\newcommand{\qgot}{\ensuremath{\mathfrak{q}}}

\newcommand{\zgot}{\ensuremath{\mathfrak{z}}}

\newcommand{\Rgot}{\ensuremath{\mathfrak{R}}}

\newcommand{\Xgot}{\ensuremath{\mathfrak{X}}}

\newcommand{\Acal}{\ensuremath{\mathcal{A}}}

\newcommand{\Ccal}{\ensuremath{\mathcal{C}}}

\newcommand{\Ecal}{\ensuremath{\mathcal{E}}}
\newcommand{\Fcal}{\ensuremath{\mathcal{F}}}

\newcommand{\Kcal}{\ensuremath{\mathcal{K}}}
\newcommand{\Lcal}{\ensuremath{\mathcal{L}}}

\newcommand{\Ocal}{\ensuremath{\mathcal{O}}}

\newcommand{\Qcal}{\ensuremath{\mathcal{Q}}}

\newcommand{\Vcal}{\ensuremath{\mathcal{V}}}

\newcommand{\Xcal}{\ensuremath{\mathcal{X}}}

\newcommand{\Qbb}{\ensuremath{\mathbb{Q}}}
\newcommand{\Nbb}{\ensuremath{\mathbb{N}}}

\newcommand{\Cbb}{\ensuremath{\mathbb{C}}}
\newcommand{\Rbb}{\ensuremath{\mathbb{R}}}
\newcommand{\Zbb}{\ensuremath{\mathbb{Z}}}

\newcommand{\chol}{\ensuremath{\Ccal_{\rm hol}}}
\newcommand{\wchol}{\ensuremath{\tilde{\Ccal}_{\rm hol}}}

\newcommand{\Ghol}{\ensuremath{\widehat{G}_{\rm hol}}}

\newcommand{\wGhol}{\ensuremath{\widehat{\tilde{G}}_{\rm hol}}}

\newcommand{\horn}{\ensuremath{\hbox{\rm Horn}}}

\newcommand{\ad}{{\ensuremath{\rm ad}}}

\newcommand{\croc}{\ensuremath{\hookrightarrow}}
\newcommand{\T}{\ensuremath{\hbox{\bf T}}}

\newcommand{\tr}{\operatorname{Tr}}

\newcommand{\Horn}{\ensuremath{{\rm Horn}_{\rm hol}}}

\def \wG {\tilde{G}}
\def \wK {\tilde{K}}
\def \wT {\tilde{T}}

\def \wg {\tilde{\ggot}}

\def \wk {\tilde{\kgot}}
\def \wt {\tilde{\tgot}}

\begin{document}

\title{Horn problem for quasi-hermitian Lie groups}

\author{Paul-Emile Paradan\footnote{IMAG, Univ Montpellier, CNRS, email : paul-emile.paradan@umontpellier.fr}}

\maketitle

\date{}


\begin{abstract}
In this paper, we prove some convexity results associated to orbit projection of non-compact real reductive Lie groups.
\end{abstract}



\tableofcontents

\section{Introduction}

This paper is concerned with convexity properties associated to orbit projection.

Let us consider two Lie groups $G\subset \wG$ with Lie algebras $\ggot\subset\wg$ and corresponding projection $\pi_{\ggot,\wg}: \wg^*\to\ggot^*$. A 
longstanding problem has been to understand how a coadjoint orbit $\tilde{\Ocal}\subset \wg^*$ decomposes under the projection 
$\pi_{\ggot,\wg}$. For this purpose, we may define 
$$
\Delta_G(\tilde{\Ocal})=\{\Ocal\in \ggot^*/G\, ; \ \Ocal\subset \pi_{\ggot,\wg}(\tilde{\Ocal})\}.
$$

When the Lie group $G$ is compact and connected, the set $\ggot^*/G$ admits a natural identification with a Weyl chamber $\tgot^*_{\geq 0}$. In this context, 
we have the well-known convexity theorem \cite{Heckman82,Atiyah82,Guillemin-Sternberg82.bis,Kirwan.84.bis,Hilgert-Neeb-Plank-94,sjamaar98,L-M-T-W}.

\begin{theo}\label{theo:convexity-1}
Suppose that $G$ is compact connected and that the projection $\pi_{\ggot,\wg}$ is proper when restricted to $\tilde{\Ocal}$. Then 
$\Delta_G(\tilde{\Ocal})=\{\xi\in \tgot^*_{\geq 0}\, ; \ G\xi\subset \pi_{\ggot,\wg}(\tilde{\Ocal})\}$ is a closed convex locally polyhedral subset of $\tgot^*$.
\end{theo}

When the Lie group $\tilde{G}$ is also compact and connected, we may consider 
\begin{equation}\label{eq:pi-G-tilde-G}
\Pi(\tilde{G},G):=\left\{(\tilde{\xi},\xi)\in \tilde{\tgot}^*_{\geq 0}\times\tgot^*_{\geq 0} ; \ G\xi\,\subset \,\pi_{\ggot,\tilde{\ggot}}\big(\tilde{G}\tilde{\xi}\,\big)\right\}.
\end{equation}

Here is another convexity theorem \cite{Horn,Klyachko,Berenstein-Sjamaar,Belkale06,BK1,Montagard-Ressayre-07,Knutson-Tao-Woodward-04,Ressayre-10}.
\begin{theo}\label{theo:convexity-2}
Suppose that $G\subset \wG$ are compact connected Lie groups. Then $\Pi(\tilde{G},G)$ is a closed convex polyhedral cone and we can parametrize its facets by cohomological means (i.e. Schubert calculus).
\end{theo} 

In this article we obtain an extension of Theorems \ref{theo:convexity-1} and \ref{theo:convexity-2} in the case where
the two groups $G$ and $\wG$ are not compact.

Let us explain in which framework we work. We suppose that $\wG/\wK$ is a Hermitian symmetric space of a non-compact type. Among the elliptic coadjoint orbits of $\wG$, some of them are naturally K\"{a}hler $\wK$-manifolds. These orbits are called the holomorphic coadjoint orbits of $\wG$. They are the strongly elliptic coadjoint orbits closely related to the holomorphic discrete series of Harish-Chandra. These orbits intersect the Weyl chamber $\tilde{\tgot}^*_{\geq 0}$ of $\wK$ into a sub-chamber $\wchol$ called the holomorphic chamber. The basic fact here is that the union 
$$
\Ccal_{\wG/\wK}^0:=\bigcup_{\tilde{a}\in\wchol} \wG\tilde{a}
$$ 
is an open invariant convex cone of $\wg^*$. See \S \ref{sec:cone-elliptic} for more details.

Let $G/K\subset \wG/\wK$ be a sub-Hermitian symmetric space of a non-compact type. As the projection $\pi_{\ggot,\wg}: \wg^*\to\ggot^*$ sends the convex cone  
$\Ccal_{\wG/\wK}^0$ inside the convex cone $\Ccal_{G/K}^0$, it is natural to study the following object reminiscent of (\ref{eq:pi-G-tilde-G}) : 
\begin{equation}
\Pi_{\rm hol}(\tilde{G},G):=\left\{(\tilde{\xi},\xi)\in \wchol\times\chol ; \ G\xi\,\subset \,\pi_{\ggot,\tilde{\ggot}}\big(\tilde{G}\tilde{\xi}\,\big)\right\}.
\end{equation}

Let $\tilde{\mu}\in \wchol$. We will also give a particular attention to the intersection of $\Pi_{\rm hol}(\tilde{G},G)$ with the hyperplane $\tilde{\xi}=\tilde{\mu}$, 
that is to say 
\begin{equation}
\Delta_G(\wG\tilde{\mu}):=\left\{\xi\in \chol ; \ G\xi\,\subset \,\pi_{\ggot,\tilde{\ggot}}\big(\tilde{G}\tilde{\mu}\,\big)\right\}.
\end{equation}

Consider the case where $G$ is embedded diagonally in $\wG:=G^s$ for $s\geq 2$. The corresponding set $\Pi_{\rm hol}(G^s,G)$ 
is called the holomorphic Horn cone, and it is defined as follows
$$
\Horn^s(G):=\Big\{(\xi_1,\cdots,\xi_{s+1})\in \chol^{s+1} ; \ G\xi_{s+1}\,\subset \sum_{j=1}^s G\xi_j\Big\}.
$$

The first result of this article is the following Theorem.

\medskip

\noindent {\bf Theorem A.}
\begin{itemize} 
\item $\Pi_{\rm hol}(\tilde{G},G)$ is a closed convex cone of $\wchol\times\chol$.
\item $\Horn^s(G)$ is a closed convex cone of $\chol^{s+1}$ for any $s\geq 2$.
\end{itemize}

\medskip

We obtain the following corollary which corresponds to a result of Weinstein \cite{Weinstein01}.

\medskip

\noindent {\bf Corollary.}{\em \ 
For any $\tilde{\mu}\in \wchol$, $\Delta_G(\wG\tilde{\mu})$ is a closed and convex subset of $\chol$.
}

\medskip

A first description of the closed convex cone $\Pi_{\rm hol}(\tilde{G},G)$ goes as follows. The quotient $\qgot$ of the tangent spaces
$\T_{e}G/K$ and $\T_{e}\wG/\wK$ has a natural structure of a Hermitian $K$-vector space. The symmetric algebra ${\rm Sym}(\qgot)$ of 
$\qgot$ defines an admissible $K$-module. The irreducible representations of $K$ (resp. $\wK$) are parametrized by a semi-group 
$\wedge^*_+:$ (resp. $\tilde{\wedge}^*_+$). For any $\lambda\in\wedge^*_+$ (resp.  $\tilde{\lambda}\in\tilde{\wedge}^*_+$), we denote by $V_\lambda^K$ (resp. $V_{\tilde{\lambda}}^{\wK}$) the irreducible representation of $K$ (resp. $\wK$) with highest weight $\lambda$ (resp. $\tilde{\lambda}$). If $E$ is a representation of 
$K$, we denote by $\left[V^K_\lambda : E\right]$ the multiplicity of $V^K_\lambda$ in $E$.

\begin{defi}\label{defi:horn-cone-general} 
\begin{enumerate}
\item $\Pi^\Zbb(\wG,G)$ is the semigroup of $\tilde{\wedge}^*_+\times\wedge^*_+$ defined by the conditions:
$$
(\tilde{\lambda},\lambda)\in \Pi^\Zbb(\wG,G)\Longleftrightarrow \Big[V^K_\lambda\,:\,V_{\tilde{\lambda}}^{\wK}\otimes {\rm Sym}(\qgot)\Big]\neq 0.
$$
\item $\Pi(\wG,G)$ is the convex cone defined as the closure of 
$\Qbb^{>0}\cdot \Pi^\Zbb(\wG,G)$.
\end{enumerate}
\end{defi}

The second result of this article is the following Theorem.

\medskip

\medskip

\noindent {\bf Theorem B.} {\em We have the equality 
\begin{equation}\label{eq:theo-B}
\Pi_{\rm hol}(\wG,G)=\Pi(\wG,G)\, \bigcap\,  \wchol\!\times\!\chol.
\end{equation}
}

\medskip

A natural question is the description of the facets of the convex cone $\Pi_{\rm hol}(\tilde{G},G)$. In order to do that, 
we consider the group $\wK$ endowed with the following $\wK\times K$-action : $(\tilde{k},k)\cdot \tilde{a}=\tilde{k}\tilde{a}k^{-1}$. 
The cotangent space $\T^*\wK$ is then a symplectic manifold equipped with a Hamiltonian action of $\wK\times K$.
We consider now the Hamiltonian $\wK\times K$-manifold $\T^*\wK\times \qgot$, and we denote by $\Delta(\T^*\wK\times \qgot)$  
 the corresponding  Kirwan polyhedron.

Let $W=N(T)/T$ be the Weyl group of $(K,T)$ and let $w_0$ be the longest Weyl group element. Define an involution $*:\tgot^*\to\tgot^*$ by 
$\xi^*=- w_0\xi$. A standard result permits to affirm that $(\tilde{\xi},\xi)\in \Pi(\wG,G)$ is an only if 
$(\tilde{\xi},\xi^*)\in\Delta(\T^*\wK\times \qgot)$  (see \S \ref{sec:affine}). 

We obtain then another version of Theorem B.

\medskip

\medskip

\noindent {\bf Theorem B, second version.} {\em An element $(\tilde{\xi},\xi)$ belongs to  $\Pi_{\rm hol}(\wG,G)$ if and only if 
$$
(\tilde{\xi},\xi)\in\wchol\times\chol\quad {\rm and}\quad (\tilde{\xi},\xi^*)\in\Delta(\T^*\wK\times \qgot).
$$
}

\medskip

Thanks to the second version of Theorem B, a natural way to describe the facets of the cone $\Pi_{\rm hol}(\wG,G)$ is to exhibit those of 
the Kirwan polyhedron $\Delta(\T^*\wK\times \qgot)$. In this later case, it can be done using Ressayre's data (see \S \ref{sec:cone:faces}).

The second version of Theorem B  permits also the following description of the convex subsets $\Delta_G(\wG\tilde{\mu})$, $\tilde{\mu}\in \wchol$. 
Let $\Delta_K(\wK\tilde{\mu}\times\overline{\qgot})$  be the Kirwan polyhedron associated to the Hamiltonian action of 
$K$ on $\wK\tilde{\mu}\times \overline{\qgot}$. Here $\overline{\qgot}$ denotes the 
$K$-module $\qgot$ with opposite complex structure.

\medskip

\medskip
 
\noindent {\bf Theorem C.} {\em For any $\tilde{\mu}\in \wchol$, we have $\Delta_G(\wG\tilde{\mu})=\Delta_K(\wK\tilde{\mu}\times\overline{\qgot})$.
}

\medskip

\medskip

Let us detail Theorem C in the case where $G$ is embedded in $\wG=G\times G$ diagonally. We denote by 
$\pgot$  the $K$-Hermitian space $\T_{e}G/K$. Let $\kappa$ be the Killing form of the Lie algebra $\ggot$. The vector space $\overline{\pgot}$ is equipped with the symplectic $2$-form $\Omega_{\bar{\pgot}}(X,Y)=-\kappa(z,[X,Y])$ and the compatible complex structure $-{\rm ad}(z)$.

Let us denote by $\Delta_K(K\mu_1\times K\mu_2\times\overline{\pgot})$ and by $\Delta_K(\overline{\pgot})$ the Kirwan polyhedrons relative to the Hamiltonian actions of 
$K$ on $K\mu_1\times K\mu_2\times\overline{\pgot}$ and on $\overline{\pgot}$. Theorem C says that for any 
$\mu_1,\mu_2\in\chol$, the convex set $\Delta_G(G\mu_1\times G\mu_2)$ is equal to the Kirwan polyhedron 
$\Delta_K(K\mu_1\times K\mu_2\times\overline{\pgot})$.

To any non-empty subset $C$ of a real vector space $E$, we may associate its asymptotic cone ${\rm As}(C) \subset E$ which is  the set formed by the limits 
$y=\lim_{k\to\infty} t_k y_k$, where $(t_k)$ is a sequence of non-negative reals converging to $0$ and $y_k\in C$.

We finally get the following description of the asymptotic cone of $\Delta_G(G\mu_1\times G\mu_2)$.

\medskip

\medskip
 
\noindent {\bf Corollary D.} {\em
For any $\mu_1,\mu_2\in\chol$, the asymptotic cone of $\Delta_G(G\mu_1\times G\mu_2)$ is equal to $\Delta_K(\overline{\pgot})$.
}

\medskip

\medskip

In \cite{pep-discrete-2008} \S 5,  we explained how to describe the cone $\Delta_K(\overline{\pgot})$ in terms of strongly orthogonal roots.

\medskip

Let us finish this introduction with few remarks on related works:

\begin{enumerate}
\item[$-$] When $G$ is compact, equal to the maximal compact subgroup $\wK$ of $\wG$, the results of Theorems B and C were already obtained by G. Deltour in his thesis 
\cite{Deltour-these,Deltour-transf-group}.

\item[$-$] In \cite{Eshmatov-Foth14}, A. Eshmatov and P. Foth proposed a description of the set $\Delta_G(G\mu_1\times G\mu_2)$. {\bf But 
their computations give not the same result as ours}. From their main result (Theorem 3.2) it follows that the asymptotic cone of $\Delta_G(G\mu_1\times G\mu_2)$ 
is equal to the intersection of the Kirwan polyhedron $\Delta_T(\overline{\pgot})$ with the Weyl chamber $\tgot^*_{\geq 0}$. But 
since $\Delta_K(\overline{\pgot})\neq \Delta_T(\overline{\pgot})\cap \tgot^*_{\geq 0}$ in general, it is in contradiction with Corollary D.
\end{enumerate}

\subsection*{Notations}

 In this paper, $G$ denotes a connected real reductive Lie group : we take here 
the convention of Knapp \cite{Knapp-book}. We have a Cartan involution $\Theta$ on $G$, such 
that the fixed point set $K:=G^\Theta$ is a connected maximal compact subgroup. We have Cartan 
decompositions at the level of Lie algebras $\ggot=\kgot\oplus \pgot$ and at the level of the group 
$G\simeq K\times \exp(\pgot)$. We denote by $b$ a $G$-invariant non-degenerate bilinear form 
on $\ggot$ that is equal to the Killing form on $[\ggot,\ggot]$, and that defines a $K$-invariant 
scalar product $(X,Y):=-b(X,\Theta(Y))$. We will use the $K$-equivariant identification 
$\xi\mapsto \tilde{\xi},\ \ggot^*\simeq \ggot$ defined by 
$(\tilde{\xi},X):=\langle\xi,X\rangle$ for $\xi\in\ggot^*$ and $X\in\ggot$.

When a Lie group $H$ acts on a manifold $N$, the stabilizer subgroup of $n\in N$ is denoted by $H_n=\{g\in G,gn=n\}$, and its Lie algebra by $\hgot_n$.
Let us define 
\begin{equation}\label{eq:dim-K}
\dim_H(\Xcal)=\min_{n\in \Xcal} \dim(\hgot_n)
\end{equation} 
for any subset $\Xcal\subset N$.

\section{The cone $\Pi_{\rm hol}(\tilde{G},G)$: first properties}\label{sec:cone}

We assume here that $G/K$ is a Hermitian symmetric space, that is to say, there exists a $G$-invariant complex structure on the manifold 
$G/K$, or equivalently there exists a $K$-invariant element $z\in\kgot$ such that $\ad(z)\vert_{\pgot}$ defines a complex structure on $\pgot$ : 
$(\ad(z)\vert_{\pgot})^2= -{\rm Id}_{\pgot}$. This condition imposes that the ranks of $G$ and $K$ are equal.

We are interested in the following closed invariant convex cone of $\ggot^*$:
$$
\Ccal_{G/K}=\left\{\xi\in\ggot^*,\langle \xi,g z\rangle \geq 0,\ \forall g\in G\right\}.
$$

\subsection{The holomorphic chamber}\label{sec:cone-elliptic}

Let $T$ be a maximal torus of $K$, with Lie algebra $\tgot$. Its dual $\tgot^*$ can be seen as the subspace of $\ggot^*$ fixed by $T$. 
Let us denote by $\ggot^*_e$ the set formed by the elliptic elements: in other words $\ggot^*_e:={\rm Ad}^*(G)\cdot \tgot^*$.  

Following \cite{Weinstein01}, we consider the invariant open subset 
$\ggot^*_{se}=\{\xi\in\ggot^*\,\vert\, G_\xi \ \mathrm{is\  compact} \}$ 
of {\em strongly elliptic} elements.  It is non-empty since the groups $G$ and $K$ have the same rank.

We start with the following basic facts.

\begin{lem}\label{lem:C-se}
\begin{itemize}
\item $\ggot^*_{se}$ is contained in $\ggot^*_{e}$.
\item  The interior $\Ccal_{G/K}^0$ of the cone $\Ccal_{G/K}$ is contained in $\ggot^*_{se}$.

\end{itemize}
\end{lem}
{\em Proof :}  The first point is due to the fact that every compact subgroup of $G$ is conjugate to a subgroup of $K$.

Let $\xi\in \Ccal_{G/K}^0$. There exists $\epsilon>0$ so that 
$$
\langle \xi +\eta,g z\rangle \geq 0,\quad \forall g\in G,\quad \forall \|\eta\|\leq \epsilon.
$$
It implies that  $|\langle \eta,g z\rangle | \leq\langle \xi, z\rangle$, $\forall g\in G_\xi$ and $\forall \|\eta\|\leq \epsilon$. In other words, the adjoint orbit 
$G_\xi\cdot z\subset \ggot$ is bounded. For any $g=e^X k$, with $(X,k)\in \pgot\times K$, a direct computation shows that $\|gz\|=\|e^X z\|\geq \|[z,X]\|=\|X\|$. Then, 
there exists $\rho>0$ such that $\|X\|\leq\rho$ if $e^X k\in G_\xi$ for some $k\in K$. It follows that the stabilizer subgroup $G_\xi$ is compact.
$\Box$

\medskip

Let $\wedge^*\subset\tgot^*$ be the weight lattice : $\alpha\in\wedge^*$ if $i\alpha$ is the 
differential of a character of $T$. Let $\Rgot\subset \wedge^*$ be the set of roots 
for the action of $T$ on $\ggot\otimes\Cbb$. We have $\Rgot=\Rgot_c\cup \Rgot_n$ 
where $\Rgot_c$ and $\Rgot_n$ are respectively the set of roots for the action of $T$ on 
$\kgot\otimes\Cbb$ and $\pgot\otimes\Cbb$. We fix a system of positive roots $\Rgot^+_c$ in 
$\Rgot_c$, and we denote by $\tgot^*_{\geq 0}$ the corresponding Weyl chamber.

We have $\pgot\otimes\Cbb=\pgot^+\oplus\pgot^-$ where the $K$-module $\pgot^\pm$ 
is equal to $\ker({\rm ad}(z)\mp i)$. 
Let $\Rgot^{\pm,z}_n$ be the set of roots for the action of $T$ 
on $\pgot^\pm$. The union 
\begin{equation}\label{eq:R-plus}
\Rgot^+=\Rgot^+_c\cup\Rgot^{+,z}_n
\end{equation} 
defines then a system of positive roots in $\Rgot$. We notice that 
$\Rgot^{+,z}_n$ is the set of roots $\beta\in \Rgot$ satisfying $\langle\beta, z\rangle=1$. Hence $\Rgot^{+,z}_n$ is invariant relatively to the action of the 
Weyl group $W=N(T)/T$.

Let us recall the following classical fact concerning the parametrization of the $G$-orbits in $\Ccal_{G/K}^0$ via the holomorphic chamber
$$
\chol:=\{\xi\in\tgot^*_{\geq 0},(\xi,\beta) >0,\ \forall \beta\in \Rgot^{+,z}_{n}\}.
$$

\medskip

\begin{prop}\label{prop:Gz}
\begin{itemize}
\item The set $\Ccal_{G/K}^0\cap \tgot^*_{\geq 0}$ is contained in $\chol$.
\item  For any compact subset $\Kcal$ of $\chol$, there exists $c_\Kcal>0$ such that 
$\langle \xi,g z\rangle \geq c_\Kcal \|g z\|$,  $\forall g\in G$, $\forall \xi\in\Kcal$.
\item The map $\Ocal\mapsto \Ocal\cap \tgot^*_{\geq 0}$ defines a bijective map between the set of $G$-orbits in $\Ccal_{G/K}^0$ and the holomorphic chamber
$\chol$.
\end{itemize}
\end{prop}
{\em Proof :} 
Let $\xi\in\tgot^*_{\geq 0}\cap\Ccal_{G/K}^0$. The first point is proved if we check that  $(\xi,\beta) >0$  for any $\beta\in\Rgot^{+,z}_n$. 
Let $X_\beta,Y_\beta\in \pgot$ such that $X_\beta+i Y_\beta\in  (\pgot\otimes\Cbb)_\beta$. We choose the following normalization~: 
the vector $h_\beta:=[X_\beta,Y_\beta]$ satisfies $\langle\beta,h_\beta\rangle=1$. We see then that $(\xi,\beta)=\frac{1}{\|h_\beta\|^2}\langle \xi,h_\beta\rangle$
for any $\xi\in \ggot^*$. Standard computation \cite{Paneitz83} gives 
$$
e^{t\,{\rm ad}(X_\beta)}z=z +({\rm cosh}(t)-1)h_\beta +{\rm sinh}(t)Y_\beta,\quad \forall t\in\Rbb.
$$
By definition, we must have $\langle \xi+\eta,e^{t\,{\rm ad}(X_\beta)} z\rangle\geq 0, \forall t\in\Rbb$, for any $\eta\in \tgot^*$ small enough. It imposes that 
$\langle \xi,h_\beta\rangle>0$. The first point is settled.

Now choose some maximal strongly orthogonal system $\Sigma\subset \Rgot^{+,z}_n$. The real span $\agot$ of the $X_\beta,\beta\in\Sigma$ is a maximal abelian subspace of $\pgot$. Hence any element $g\in G$ can be written $g=k e^{X(t)} k'$ with $X(t)=\sum_{\beta\in\Sigma}t_\beta X_\beta$ and $k,k'\in K$. We get
\begin{equation}\label{eq:Gz}
gz=k\left(z +\sum_{\beta\in\Sigma}({\rm cosh}(t_\beta)-1)h_\beta +\sum_{\beta\in\Sigma}{\rm sinh}(t_\beta)Y_\beta\right)
\end{equation}
and 
$$
\langle \xi,g z\rangle=\langle k^{-1}\xi,z\rangle + \sum_{\beta\in\Sigma}({\rm cosh}(t_\beta)-1) \langle k^{-1}\xi,h_\beta\rangle.
$$

For any $\xi\in \chol$, we define $c_\xi := \min_{\beta\in \Rgot^{+,z}_n}\langle \xi,h_\beta\rangle>0$. Let $\pi:\kgot^*\to \tgot^*$ be the projection. We have 
$\langle k^{-1}\xi,z\rangle=\langle \pi(k^{-1}\xi),z\rangle$ and $\langle k^{-1}\xi,h_\beta\rangle=\langle \pi(k^{-1}\xi),h_\beta\rangle$. 
The convexity theorem of Kostant tell us that $\pi(k^{-1}\xi)$ belongs to the convex hull of $\{w\xi, w\in W\}$. It follows that 
$\langle k^{-1}\xi,z\rangle\geq \langle \xi,z\rangle>0 \quad {\rm and}\quad \langle k^{-1}\xi,h_\beta\rangle\geq c_\xi>0$. 
We obtain then that $\langle \xi,g z\rangle \geq \frac{1}{2}\min (\langle \xi,z\rangle, c_\xi) e^{\|X(t)\|}$ for any $\xi\in\chol$, where $\|X(t)\|=\sup |t_\beta|$. 
From (\ref{eq:Gz}), it is not difficult to see that there exists $C>0$ such that $\|gz\|\leq C e^{\|X(t)\|}$ for any $g=k e^{X(t)} k'$. 

Let $\Kcal$ be a compact subset of $\chol$. Take $c_\Kcal=\frac{1}{2C}\min (\min_{\xi\in\Kcal}\langle \xi,z\rangle, \min_{\xi\in\Kcal}c_\xi)>0$. The previous computations show that $\langle \xi,g z\rangle \geq c_\Kcal \|g z\|,\quad \forall g\in G,\ \forall \xi\in\Kcal$. The second point is proved.

The first two points show that $\tgot^*_{\geq 0}\cap\Ccal_{G/K}^0=\chol$. The third point follows then from the fact that $\Ccal_{G/K}^0$ 
is contained in $\ggot^*_{e}$ (see Lemma \ref{lem:C-se}).
$\Box$

\medskip

We have a canonical map ${\rm p}: \Ccal_{G/K}^0\to \chol$ defined by the relations $G\xi\cap \tgot^*_{\geq 0}=\{ {\rm p}(\xi)\}$, $\forall \xi\in \Ccal_{G/K}^0$. The following result is needed in \S \ref{sec:proof-theo-A}.
\begin{lem}\label{lem:p-continu}
The map ${\rm p}$ is continuous.
\end{lem}
{\em Proof :} Let $(\xi_n)$ be a sequence of $\Ccal_{G/K}^0$ converging to $\xi_\infty\in \Ccal_{G/K}^0$. Let $\xi'_n ={\rm p}(\xi_n)$ and 
$\xi'_\infty ={\rm p}(\xi_\infty)$: we have to prove that the sequence $(\xi'_n)$ converges to $\xi'_\infty$. We choose elements $g_n,g_\infty\in G$ such that 
$\xi_n=g_n\xi'_n,\forall n$, and $\xi_\infty=g_\infty\xi'_\infty$.

First we notice that $-b(\xi_n,\xi_n)=\|\xi'_n\|^2$, hence the sequence $(\xi'_n)$ is bounded. We will now prove that the sequence $(g_n)$ is bounded.

 Let $\epsilon>0$ such that 
$\langle\xi_\infty+\eta,gz\rangle \geq 0$, $\forall g\in G$, $\forall \|\eta\|\leq \epsilon$. If $\|\xi-\xi_\infty\|\leq \epsilon/2$, we write 
$\xi=\frac{1}{2}(\xi_\infty+ 2(\xi-\xi_\infty))+\frac{1}{2}\xi_\infty$, and then 
$$
\langle\xi,gz\rangle= \frac{1}{2}\langle\xi_\infty+ 2(\xi-\xi_\infty),gz\rangle +\frac{1}{2}\langle\xi_\infty,gz\rangle\geq \frac{1}{2}\langle\xi_\infty,gz\rangle,
\quad \forall g\in G.
$$
Now we have $\langle\xi'_n,z\rangle=\langle\xi_n,g_nz\rangle\geq \frac{1}{2}\langle\xi_\infty,g_nz\rangle$ if $n$ is large enough. This shows that the sequence 
$\langle\xi_\infty,g_nz\rangle$ is bounded. If we use the second point of Proposition \ref{prop:Gz},
we can conclude that the sequence $(g_n)$ is bounded.

Let $(\xi'_{\phi(n)})$ be a subsequence converging to $\ell\in\tgot^*_{\geq 0}$. Since  $(g_{\phi(n)})$ is bounded, there exists a 
subsequence $(g_{\phi\circ\psi(n)})$ converging to $h\in G$. From the relations $\xi_{\phi\circ\psi(n)}=$ \break 
$g_{\phi\circ\psi(n)}\xi'_{\phi\circ\psi(n)},\forall n\in\Nbb$, we obtain $\xi_\infty=h\ell$. Then $\ell={\rm p}(\xi_\infty)=
\xi_\infty'$. Since every subsequence of $(\xi'_{n})$ has a subsequential limit to $\xi_\infty'$, then the sequence $(\xi'_n)$ converges to $\xi'_\infty$.
$\Box$

\subsection{The cone $\Pi_{\rm hol}(\tilde{G},G)$ is closed}\label{sec:cone-closed}

We suppose that $G/K$ is a complex submanifold of a Hermitian symmetric space $\wG/\wK$. In other words,  $\tilde{G}$ is a reductive real Lie group  
such that $G\subset \tilde{G}$ is a closed subgroup preserved by the Cartan involution, and $\tilde{K}$ is a maximal compact subgroup of $\tilde{G}$ 
containing $K$. We denote by $\tilde{\ggot}$ and $\tilde{\kgot}$ the Lie algebras of $\tilde{G}$ and $\tilde{K}$ respectively. We suppose that there exists a $\tilde{K}$-invariant element $z\in\kgot$ such that $\ad(z)\vert_{\tilde{\pgot}}$ defines a complex structure on $\tilde{\pgot}$ : $(\ad(z)\vert_{\tilde{\pgot}})^2= -Id_{\tilde{\pgot}}$.

Let $\Ccal_{\wG/\wK}\subset \wg^*$ be the closed invariant cone associated to the Hermitian symmetric space $\wG/\wK$. We start with the following key fact.
\begin{lem}
The projection $\pi_{\ggot,\wg}: \wg^*\to\ggot^*$ sends $\Ccal_{\wG/\wK}^0$ into $\Ccal_{G/K}^0$.
\end{lem}
{\em Proof :} 
Let $\tilde{\xi}\in \Ccal_{\wG/\wK}^0$ and $\xi=\pi_{\ggot,\wg}(\tilde{\xi})$. Then $\langle \tilde{\xi}+\tilde{\eta},\tilde{g} z\rangle \geq 0,\ \forall \tilde{g}\in \wG$ 
if $\tilde{\eta}\in\wg^*$ is small enough. It follows that $\langle \xi+\pi_{\ggot,\wg}(\tilde{\eta}),g z\rangle=\langle \tilde{\xi}+\tilde{\eta},g z\rangle \geq 0,\ \forall g\in G$ 
 if $\tilde{\eta}$ is small enough. Since $\pi_{\ggot,\wg}$ is an open map, we can conclude that $\xi\in\Ccal_{G/K}^0$.
$\Box$

\medskip

Let $\wT$ be a maximal torus of $\wK$, with Lie algebra $\wt$. The $\wG$-orbits in the interior of $\Ccal_{\wG/\wK}$ are parametrized by the holomorphic chamber 
$\wchol\subset\wt^*$.
The previous Lemma says that the projection $\pi_{\ggot,\wg}(\tilde{\Ocal})$ of any $\wG$-orbit $\tilde{\Ocal}\subset \Ccal_{\wG/\wK}^0$ is the union of $G$-orbits
$\Ocal\subset \Ccal^0_{G/K}$.  So it is natural to study the following object
\begin{equation}
\Pi_{\rm hol}(\tilde{G},G):=\left\{(\tilde{\xi},\xi)\in \wchol\times\chol ; \ G\xi\,\subset \,\pi_{\ggot,\tilde{\ggot}}\big(\tilde{G}\tilde{\xi}\,\big)\right\}.
\end{equation}

Here is a first result. 

\begin{prop}\label{prop:pi-z-compact}
$\Pi_{\rm hol}(\tilde{G},G)$ is a closed cone of $\wchol\times\chol$.
\end{prop}

{\em Proof :} 
Suppose that a sequence $(\tilde{\xi}_n,\xi_n)\in \Pi_{\rm hol}(\tilde{G},G)$ converge to $(\tilde{\xi}_\infty,\xi_\infty)\in\wchol\times\chol$.
By definition, there exists a sequence $(\tilde{g}_n,g_n)\in \tilde{G}\times G$ such that 
$g_n \xi_n=\pi_{\ggot,\tilde{\ggot}}(\tilde{g}_n\tilde{\xi}_n)$. Let $\tilde{h}_n:=g_n^{-1}\tilde{g}_n$, so that $\xi_n=\pi_{\ggot,\tilde{\ggot}}(\tilde{h}_n\tilde{\xi}_n)$ and  
$\langle \tilde{h}_n\tilde{\xi}_n, z\rangle=\langle \xi_n, z\rangle$. We use now that the sequence $\langle \xi_n, z\rangle$ is bounded, and that the sequence 
$\tilde{\xi}_n$ belongs to a compact subset of $\wchol$. Thanks to the second point of Proposition \ref{prop:Gz}, these facts imply that $\|\tilde{h}_n^{-1}z\|$ is a bounded sequence. Hence $\tilde{h}_n$ admits a subsequence converging to $\tilde{h}_\infty$. So we get $\xi_\infty=\pi_{\ggot,\tilde{\ggot}}(\tilde{h}_\infty\tilde{\xi}_\infty)$, and that proves that $(\tilde{\xi}_\infty,\xi_\infty)\in\Pi_{\rm hol}(\tilde{G},G)$.
$\Box$

\subsection{Rational and weakly regular points}\label{sec:rational}

Let $(M,\Omega)$ be a symplectic manifold. We suppose that there exists a line bundle $\Lcal$ with connection $\nabla$ that prequantizes the $2$-form $\Omega$: in other words, $\nabla^2=- i\, \Omega$. Let $K$ be a compact connected Lie group acting on $\Lcal\to M$, and leaving the connection  invariant. Let $\Phi_K:M\to\kgot^*$ be the moment map defined by Kostant's relations
\begin{equation}\label{eq:kostant}
L_X-\nabla_X=i\langle \Phi_K,X\rangle,\quad \forall X\in\kgot.
\end{equation}
Here $L_X$ is the Lie derivative acting on the sections of $\Lcal\to M$.

Remark that relations (\ref{eq:kostant}) imply, via the equivariant Bianchi formula, the relations 
\begin{equation}\label{eq:symplectic-moment-map}
\iota(X_M)\Omega=- d\langle \Phi_K,X\rangle,\quad \forall X\in\kgot,
\end{equation}
where $X_M(m) :=\frac{d}{dt}\vert_{t=0}e^{-tX}m$ is the vector field on $M$ generated by $X\in\kgot$.

\begin{defi}\label{defi:weakly-regular}
Let $\dim_K(M):=\min_{m\in M}\dim\kgot_m$. An element $\xi\in\kgot^*$ is a weakly-regular value of $\Phi_K$ if for all $m\in\Phi_K^{-1}(\xi)$ we have 
$\dim\kgot_m=\dim_K(M)$.
\end{defi}

When $\xi\in\kgot^*$ is a weakly-regular value of $\Phi_K$, the constant rank theorem tells us that $\Phi_K^{-1}(\xi)$ is a submanifold of $M$ stable under the action of the stabilizer subgroup 
$K_\xi$. We see then that the reduced space
\begin{equation}\label{eq:symplectic-reduction}
M_\xi:=\Phi_K^{-1}(\xi)/K_\xi
\end{equation}
is a symplectic orbifold.

Let $T\subset K$ be a maximal torus with Lie algebra $\tgot$. We consider the lattice $\wedge:=\frac{1}{2\pi}\ker(\exp:\tgot\to T)$ and the dual lattice 
$\wedge^*\subset \tgot^*$ defined by $\wedge^*=\hom(\wedge,\Zbb)$. We remark that $i\eta$ is a differential of a character of $T$ if and only if $\eta\in\wedge^*$. 
The $\Qbb$-vector space generated by the lattice $\wedge^*$ is denoted by $\tgot^*_\Qbb$: the vectors belonging to 
$\tgot^*_\Qbb$ are designed as rational. An affine subspace $V\subset \tgot^*$ is called rational if it is affinely generated by its rational points.

We also fix a closed positive Weyl chamber $\tgot^*_{\geq 0}$. For each relatively open face $\sigma\subset \tgot^*_{\geq 0}$, 
the stabilizer $K_\xi$ of points $\xi\in\sigma$ under the coadjoint action, does not depend on $\xi$, and will be denoted $K_\sigma$. 
The Lie algebra $\ggot_\sigma$ decomposes into its semi-simple and central parts $\kgot_\sigma=[\kgot_\sigma,\kgot_\sigma]\oplus\zgot_\sigma$. 
The subspace $\zgot_\sigma^*$ is defined to be the annihilator of $[\kgot_\sigma,\kgot_\sigma]$, or equivalently the fixed point set of the coadjoint 
$K_\sigma$ action. Notice that $\zgot_\sigma^*$ is a rational subspace of  $\tgot^*$, and that the face $\sigma$ is an open cone of $\zgot_\sigma^*$,
 
We suppose that the moment map $\Phi_K$ is {\em proper}. The Convexity Theorem \cite{Atiyah82,Guillemin-Sternberg82.bis,Kirwan.84.bis,sjamaar98,L-M-T-W} tells 
us that  $\Delta_K(M):={\rm Image}(\Phi_K)\, \bigcap\, \tgot_{\geq 0}^*$ is a closed, convex, locally polyhedral set.

\begin{defi}
We denote by $\Delta_K(M)^0$ the subset of $\Delta_K(M)$ formed by the {\em weakly-regular values} of the moment map 
$\Phi_K$ contained in $\Delta_K(M)$.
\end{defi}

We will use the following remark in the next sections.

\begin{lem}\label{lem:rational-dense}
The subset $\Delta_K(M)^0\cap \tgot_\Qbb^*$ is dense in $\Delta_K(M)$.
\end{lem}
{\em Proof :} 
Let us first explain why $\Delta_K(M)^0$ is a dense open subset of $\Delta_K(M)$. There exists a unique open face $\tau$ of the Weyl chamber $\tgot^*_{\geq 0}$ such as
$\Delta_K(M)\cap\tau$ is dense in $\Delta_K(M)$~: $\tau$ is called the {\em principal} face in \cite{L-M-T-W}. 
The Principal-cross-section Theorem \cite{L-M-T-W} tells us that $Y_\tau :=\Phi^{-1}(\tau)$ is a symplectic  $K_\tau$-manifold, with a trivial action of 
$[K_\tau,K_\tau]$. The line bundle $\Lcal_\tau:=\Lcal\vert_{Y_\tau}$ prequantizes the symplectic structure on $Y_\tau$, and relations 
(\ref{eq:symplectic-moment-map}) show that $[K_\tau,K_\tau]$ acts trivially on $\Lcal_\tau$. Moreover the restriction of $\Phi_K$ on 
$Y_\tau$ is the moment map $\Phi_\tau: Y_\tau\to \zgot_\tau^*$ associated to the action of the torus $Z_\tau=\exp(\zgot_\tau)$ on $\Lcal_\tau$.

Let $I\subset \zgot_\tau^*$ be the smallest affine subspace containing $\Delta_K(M)$. Let $\zgot_I\subset\zgot_\tau$ be the annihilator of 
the subspace parallel to $I$ : relations (\ref{eq:symplectic-moment-map}) show that $\zgot_I$ is the generic infinitesimal stabilizer 
of the $\zgot_\tau$-action  on $Y_\tau$. Hence $\zgot_I$ is the Lie algebra of the torus $Z_I:=\exp(\zgot_I)$.

We see then that any regular value of $\Phi_\tau: Y_\tau\to I$, viewed as a map with codomain $I$, is a weakly-regular value of the moment map $\Phi_K$. It explains  
why $\Delta_K(M)^0$ is a dense open subset of $\Delta_K(M)$.

The convex set $\Delta_K(M)\cap\tau$ is equal to $\Delta_{Z_\tau}(Y_\tau):={\rm Image}(\Phi_\tau)$, it is sufficient to check that 
$\Delta_{Z_\tau}(Y_\tau)^0\cap \tgot_\Qbb^*$ is dense in  $\Delta_{Z_\tau}(Y_\tau)$. The subtorus $Z_I\subset Z_\tau$ acts trivially on $Y_\tau$, and it acts on the line bundle $\Lcal_\tau$ 
through a character $\chi$. Let $\eta\in \wedge^*\cap\tgot^*_\tau$ such that $d\chi=i\eta\vert_{\zgot_I}$. The affine subspace $I$ which is equal to $\eta+(\zgot_I)^{\perp}$ is rational. 
Since the open subset $\Delta_{Z_\tau}(Y_\tau)^0$ generates the rational affine subspace $I$, we can conclude that $\Delta_{Z_\tau}(Y_\tau)^0\cap \tgot_\Qbb^*$ 
is dense in  $\Delta_{Z_\tau}(Y_\tau)$.
$\Box$

\subsection{Weinstein's theorem}\label{sec:weinstein}

Let $\tilde{a}\in \wchol$. Consider the Hamiltonian action of the group $G$ on the coadjoint orbit $\wG\tilde{a}$. The moment map 
$\Phi_G^{\tilde{a}}:\wG\tilde{a}\to\ggot^*$ corresponds to the restriction of the projection $\pi_{\ggot,\tilde{\ggot}}$ to $\wG\tilde{a}$. In this setting,  
the following conditions holds :
\begin{enumerate}
\item The action of $G$ on $\wG\tilde{a}$ is proper.
\item  The moment map $\Phi_G^{\tilde{a}}$ is a proper map since the map $\langle\Phi_G^{\tilde{a}},z\rangle$ is proper (see the second point of Proposition 
\ref{prop:Gz}).
\end{enumerate}

Conditions 1. and 2. imposes that the image of $\Phi_G^{\tilde{a}}$ is contained in the open subset 
$\ggot^*_{se}$ of strongly elliptic elements \cite{pep-jems}. Thus, the $G$-orbits contained in the image of $\Phi_G^{\tilde{a}}$ are parametrized by 
the following subset of the holomorphic chamber $\chol$ :
$$
\Delta_G(\wG\tilde{a}):={\rm Image}(\Phi_G^{\tilde{a}})\, \bigcap \, \tgot^*_{\geq 0}.
$$
We notice that $\Pi_{\rm hol}(\tilde{G},G)=\bigcup_{\tilde{a}\in \wchol}\{\tilde{a}\}\times\Delta_G(\wG\tilde{a})$.

Like in Definition \ref{defi:weakly-regular}, an element $\xi\in\ggot^*$ is a {\em weakly-regular} value of $\Phi_G^{\tilde{a}}$ if for all 
$m\in(\Phi_G^{\tilde{a}})^{-1}(\xi)$ we have 
$\dim\ggot_m=\min_{x\in \wG\tilde{a}}\dim(\ggot_x)$. We denote by  $\Delta_G(\wG\tilde{a})^0$ the set of elements $\xi\in\Delta_G(\wG\tilde{a})$ 
that are weakly-regular for $\Phi_G^{\tilde{a}}$.

\begin{theo}[Weinstein]\label{theo:weinstein}
For any $\tilde{a}\in \wchol$, $\Delta_G(\wG\tilde{a})$ is a closed convex subset contained in $\chol$. 
\end{theo}
{\em Proof :} 
We recall briefly the arguments of the proof (see \cite{Weinstein01} or \cite{pep-jems}[\S 2]). Under Conditions 1. and 2., one checks easily that $Y_{\tilde{a}}:=(\Phi_G^{\tilde{a}})^{-1}(\kgot^*)$ is a smooth 
$K$-invariant symplectic submanifold of $\wG\tilde{a}$ such that 
\begin{equation}\label{eq:slice-Y}
\wG\tilde{a}\simeq G\times_K Y_{\tilde{a}}.
\end{equation}
The moment map $\Phi_K^{\tilde{a}}:Y_{\tilde{a}}\to\kgot^*$, which corresponds to the restriction of  the map $\Phi_G^{\tilde{a}}$ to $Y_{\tilde{a}}$, 
is a proper map. Hence the Convexity Theorem tells us that $\Delta_K(Y_{\tilde{a}}):={\rm Image}(\Phi_K^{\tilde{a}})\bigcap \tgot^*_{\geq 0}$ is a closed, convex, locally polyhedral set. Thanks to the isomorphism (\ref{eq:slice-Y}) we see that $\Delta_G(\wG\tilde{a})$ coincides with the 
closed convex subset $\Delta_K(Y_{\tilde{a}})$. The proof is completed.
$\Box$

\medskip

The next Lemma is used in \S \ref{sec:proof-theo-A}.

\begin{lem}\label{lem:rational-dense-coadjoint}
Let $\tilde{a}\in \wchol$ be a rational element. Then $\Delta_G(\wG\tilde{a})^0\cap \tgot_\Qbb^*$ is dense in $\Delta_G(\wG\tilde{a})$.
\end{lem}
{\em Proof :} 
Thanks to (\ref{eq:slice-Y}), we know that $\Delta_G(\wG\tilde{a})=\Delta_K(Y_{\tilde{a}})$. Relation (\ref{eq:slice-Y}) shows also that 
$\Delta_G(\wG\tilde{a})^0=\Delta_K(Y_{\tilde{a}})^0$. Let $N\geq 1$ such that $\tilde{\mu}=N\tilde{a}\in \wedge^*\cap \chol$. The stabilizer subgroup $\wG_{\tilde{\mu}}$ is compact, equal to $\wK_{\tilde{\mu}}$. Let us denote by $\Cbb_{\tilde{\mu}}$ the one-dimensional representation of $\wK_{\tilde{\mu}}$ associated to  $\tilde{\mu}$. 
The convex set $\Delta_G(\wG\tilde{a})$ is equal to $\frac{1}{N}\Delta_G(\wG\tilde{\mu})$, so it is sufficient to check that 
$\Delta_G(\wG\tilde{\mu})^0\cap \tgot_\Qbb^*=\Delta_K(Y_{\tilde{\mu}})^0\cap \tgot_\Qbb^*$ is dense in $\Delta_G(\wG\tilde{\mu})=\Delta_K(Y_{\tilde{\mu}})$. The coadjoint orbit $\wG\tilde{\mu}$ is prequantized by the line bundle $\wG\times_{K_{\tilde{\mu}}}\Cbb_{\tilde{\mu}}$ and the symplectic slice $Y_{\tilde{\mu}}$ is prequantized by the line bundle $\Lcal_{\tilde{\mu}}:=\wG\times_{K_{\tilde{\mu}}}\Cbb_{\tilde{\mu}}\vert_{Y_{\tilde{\mu}}}$.
Thanks to Lemma \ref{lem:rational-dense}, we know that $\Delta_K(Y_{\tilde{\mu}})^0\cap \tgot_\Qbb^*$ is dense in $\Delta_K(Y_{\tilde{\mu}})$ : the 
proof is complete. $\Box$

\section{Holomorphic discrete series}\label{sec:holo-series}

\subsection{Definition}\label{sec:holo-series-def}

We return to the framework of \S \ref{sec:cone-elliptic}. We recall the notion of holomorphic discrete series representations associated to a 
Hermitian symmetric spaces $G/K$. Let us introduce
$$
\chol^\rho:=\left\{ \xi\in\tgot^*_{\geq 0} \vert\ (\xi,\beta)\geq (2\rho_n,\beta) , \  \forall \beta\in \Rgot_n^{+,z}\right\}
$$
where  $2\rho_n=\sum_{\beta\in\Rgot_n^{+,z}}\beta$ is $W$-invariant. 

\begin{lem} 
\begin{enumerate}
\item We have $\chol^\rho\subset\chol$.
\item For any $\xi\in\chol$ there exists $N\geq 1$ such that $N\xi\in\chol^\rho$.
\end{enumerate}
\end{lem}
{\em Proof :} 
The first point is due to the fact that $(\beta_0,\beta_1)\geq 0$ for any $\beta_0,\beta_1\in \Rgot_n^{+,z}$. The second point is obvious.
$\Box$

We will be interested in the following subset of dominants weights :
$$
\Ghol:=\wedge^*_+\bigcap \chol^\rho.
$$
Let ${\rm Sym}(\pgot)$ be the symmetric algebra of the complex $K$-module $(\pgot,\mathrm{ad}(z))$.

\begin{theo}[Harish-Chandra]
For any $\lambda\in\Ghol$, there exists an irreducible unitary representation of $G$, denoted $V^G_\lambda$, such that 
the vector space of $K$-finite vectors is $V^G_\lambda\vert_K:=V_\lambda^K\otimes {\rm Sym}(\pgot)$.  
\end{theo}

The set $V^G_\lambda,\lambda\in\Ghol$ corresponds to the holomorphic discrete series representations associated to the complex structure ${\rm ad}(z)$.

\subsection{Restriction}\label{sec:holo-series-restriction}

We come back to the framework of \S \ref{sec:cone-closed}. We consider two compatible Hermitian symmetric spaces $G/K\subset \wG/\wK$, and we look after the restriction of holomorphic discrete series representations of $\wG$ to the subgroup $G$. 

Let $\tilde{\lambda}\in \wGhol$. Since the representation $V^{\wG}_{\tilde{\lambda}}$ is discretely admissible relatively the circle group $\exp(\Rbb z)$, it is also discretely admissible relatively to $G$. We can be more precise  \cite{Jakobsen-Vergne, Martens,Toshi-JFA98} :

\begin{prop}\label{prop:restriction-V-lambda-G-prime}
We have an Hilbertian direct sum 
$$
V^{\wG}_{\tilde{\lambda}}\vert_{G}=\bigoplus_{\lambda\in \Ghol} m_{\tilde{\lambda}}^{\lambda} \ V^{G}_\lambda,
$$
where the multiplicity $m_{\tilde{\lambda}}^{\lambda} :=[V^{G}_\lambda:V^{\wG}_{\tilde{\lambda}}]$ is finite for any $\lambda$.
\end{prop}

The Hermitian $\wK$-vector space $\tilde{\pgot}$, when restricted to the $K$-action, admits an orthogonal decomposition $\tilde{\pgot}=\pgot\oplus \qgot$. Notice that the symmetric algebra ${\rm Sym}(\qgot)$ is an admissible $K$-module.

In \cite{Jakobsen-Vergne}, H.P. Jakobsen and M. Vergne obtained the following nice characterization of the multiplicities $[V^{G}_\lambda:V^{\wG}_{\tilde{\lambda}}]$.
Another proof is given in \cite{pep-jems}, \S 4.4.
\begin{theo}[Jakobsen-Vergne]\label{theo:JV}
Let $(\tilde{\lambda},\lambda)\in \wGhol\times\Ghol$. The multiplicity $[V^{G}_\lambda:V^{\wG}_{\tilde{\lambda}}]$ is equal to the multiplicity of the 
representation $V^{K}_\lambda$ in ${\rm Sym}(\qgot)\otimes V^{\wK}_{\tilde{\lambda}}\vert_{K} $. 
\end{theo}

\subsection{Discrete analogues of $\Pi_{\rm hol}(\tilde{G},G)$ }\label{sec:LR}

We define the following discrete analogues of the cone $\Pi_{\rm hol}(\tilde{G},G)$:

\begin{equation}
\Pi_{\rm hol}^\Zbb(\tilde{G},G):=\left\{(\tilde{\lambda},\lambda)\in \wGhol\times\Ghol\ ; \ [V^{G}_\lambda:V^{\wG}_{\tilde{\lambda}}]\neq 0 \right\}.
\end{equation}
and 
\begin{equation}
\Pi_{\rm hol}^\Qbb(\tilde{G},G):=\left\{(\tilde{\xi},\xi)\in \wchol\times\chol\ ; \ \exists N\geq 1,\ (N\xi,N\tilde{\xi})\in \Pi_{\rm hol}^\Zbb(\tilde{G},G)\right\}.
\end{equation}

We have the following key fact.

\begin{prop}\label{prop:monoid}
${}$
\begin{itemize}
\item $\Pi_{\rm hol}^\Zbb(\tilde{G},G)$ is a subset of $ \tilde{\wedge}^* \times \wedge^*$ stable under the addition. 
\item $\Pi_{\rm hol}^\Qbb(\tilde{G},G)$ is a $\Qbb$-convex cone of the $\Qbb$-vector space $\tilde{\tgot}^*_\Qbb \times \tgot^*_\Qbb$.
\end{itemize}
\end{prop}
{\em Proof :} 
Suppose that $a_1:=(\tilde{\lambda}_1,\lambda_1)$ and $a_2:=(\tilde{\lambda}_2, \lambda_2)$ belongs to 
$\Pi_{\rm hol}^\Zbb(\tilde{G},G)$. Thanks to Theorem \ref{theo:JV}, we know that the $K$-modules ${\rm Sym}(\qgot)\otimes (V^K_{\lambda_j})^*\otimes V^{\wK}_{\tilde{\lambda}_j}\vert_{K}$ possess a non-zero invariant vector $\phi_j$, for $j=1,2$. 

Let $\mathbb{X}:=\overline{K/T}\times \wK/\wT$ be the product of flag manifolds. The complex structure is normalized so that 
$\T _{([e],[\tilde{e}])}\mathbb{X}\simeq \ngot_-\oplus \tilde{\ngot}_+$, where $\ngot_-=\sum_{\alpha<0}(\kgot_\Cbb)_\alpha$ and 
$\tilde{\ngot}_+=\sum_{\tilde{\alpha}>0}(\tilde{\kgot}_\Cbb)_{\tilde{\alpha}}$. 
We associate to each data $a_j$, the holomorphic line bundle $\Lcal_j:=K\times_T \Cbb_{-\lambda_j}\boxtimes \wK\times_{\wT} \Cbb_{-\tilde{\lambda}_j}$ on 
$\mathbb{X}$. Let $H^0(\mathbb{X},\Lcal_j)$ be the space of holomorphic sections of the line bundle $\Lcal_j$. The Borel-Weil Theorem tell us that 
$H^0(\mathbb{X},\Lcal_j)\simeq (V^K_{\lambda_j})^*\otimes V^{\wK}_{\tilde{\lambda}_j}\vert_{K}$, $\forall j\in\{1,2\}$.

We have $\phi_j\in \left[{\rm Sym}(\qgot)\otimes H^0(\mathbb{X},\Lcal_j)\right]^K$, $\forall j$, and then $\phi_1\phi_2\in {\rm Sym}(\qgot)\otimes H^0(\mathbb{X},\Lcal_1\otimes \Lcal_2)$ is a non-zero invariant vector. Hence 
$[{\rm Sym}(\qgot)\otimes (V^K_{\lambda_1+\lambda_2})^*\otimes V^{\wK}_{\tilde{\lambda}_1+\tilde{\lambda}_2}\vert_{K}]^K\neq 0$.
Thanks to Theorem \ref{theo:JV}, we can conclude that $a_1+a_2=(\tilde{\lambda}_1+\tilde{\lambda}_2,\lambda_1+\lambda_2)$ belongs to 
$\Pi_{\rm hol}^\Zbb(\tilde{G},G)$. The first point is proved. From the first point, one checks easily that 
\begin{itemize}
\item[-] $\Pi_{\rm hol}^\Qbb(\tilde{G},G)$ is stable under addition,
\item[-] $\Pi_{\rm hol}^\Qbb(\tilde{G},G)$ is stable by expansion by a non-negative rational number.
\end{itemize}
The second point is settled.
$\Box$

\subsection{Riemann-Roch numbers}\label{sec:RR-numbers}
 We come back to the framework of \S \ref{sec:rational}.

 We associate to a dominant weight $\mu\in\wedge^*_{+}$ the (possibly singular) symplectic reduced space 
 $M_\mu:=\Phi_K^{-1}(\mu)/K_\mu$, and the (possibly singular) line bundle 
 over $M_\mu$~:
 $$
 \Lcal_{\mu}:=\left(\Lcal\vert_{\Phi_K^{-1}(\mu)}\otimes \Cbb_{-\mu}\right)/K_\mu.
 $$

 \underline{Suppose first that $\mu$ is a weakly-regular value of $\Phi_K$.} Then $M_\mu$ is an orbifold equipped with a symplectic structure 
 $\Omega_\mu$, and $\Lcal_{\mu}$ is a line orbi-bundle over $M_\mu$ that prequantizes the symplectic structure. 
 By choosing an almost complex structure on $M_\mu$ compatible with $\Omega_\mu$, we get a decomposition 
 $\wedge \T^* M_\mu\otimes \Cbb = \oplus_{i,j}\wedge^{i,j} \T^* M_\mu$ of the bundle of differential forms. 
 Using Hermitian structure in the tangent bundle $\T M_\mu$ of $M_\mu$, and in the fibers of $\Lcal_{\mu}$, we define a 
Dolbeaut-Dirac operator
$$
 D^+_{\mu} : \Acal^{0,+}(M_\mu,\Lcal_{\mu})\longrightarrow \Acal^{0,-}(M_\mu,\Lcal_{\mu})
$$
where $\Acal^{i,j}(M_\mu,\Lcal_{\mu})=\Gamma(M_\mu,\wedge^{i,j} \T^* M_\mu\otimes \Lcal_{\mu})$.
 
\begin{defi}
Let $\mu\in\wedge^*_{+}$ be a weakly-regular value of the moment map $\Phi_K$. The Riemann-Roch number 
$RR(M_\mu,\Lcal_{\mu})\in\Zbb$ is defined as the index of the elliptic operator $D^+_{\mu}$:
$RR(M_\mu,\Lcal_{\mu})= \dim(ker(D^+_{\mu})) -\dim(coker(D^+_{\mu}))$.
\end{defi}

\medskip

\underline{Suppose that $\mu\notin \Delta_K(M)$.} Then $M_\mu=\emptyset$ and we take $RR(M_\mu,\Lcal_{\mu})=0$.

\medskip
 
\underline{Suppose now that $\mu\in \Delta_K(M)$ is not (necessarily) a weakly-regular value of $\Phi_K$.} 
Take a small element $\epsilon\in\tgot^*$ such that $\mu+\epsilon$ is a weakly-regular value of $\Phi_K$ belonging to $\Delta_K(M)$. 
We consider the symplectic orbifold $M_{\mu+\epsilon}$: if $\epsilon$ is small enough,  
$$
\Lcal_{\mu,\epsilon}:=\left(\Lcal\vert_{\Phi_K^{-1}(\mu+\epsilon)}\otimes \Cbb_{-\mu}\right)/K_{\mu+\epsilon}.
$$ 
is a line orbi-bundle over $M_{\mu+\epsilon}$ .

We have the following important result (see \S 3.4.3 in \cite{pep-vergne-ams}).
\begin{prop}Let $\mu\in \Delta_K(M)\cap \wedge^*$. 
The Riemann-Roch number  $RR(M_{\mu+\epsilon},\Lcal_{\mu,\epsilon})$ do not depend on the choice of $\epsilon$ small enough so that 
$\mu+\epsilon\in \Delta_K(M)$ is a weakly-regular value of $\Phi_K$.
\end{prop}

We can now introduce the following definition.

\begin{defi}\label{def:RR-numbers}
Let $\mu\in\wedge^*_{+}$.  We define 
$$
\Qcal(M_\mu,\Omega_\mu)=
\begin{cases}
0 \hspace{29mm} {\rm if }\quad  \mu\notin \Delta_K(M),\\
RR(M_{\mu+\epsilon},\Lcal_{\mu,\epsilon})\hspace{4mm}  {\rm if }\quad \mu\in \Delta_K(M).
\end{cases}
$$
Above, $\epsilon$ is chosen small enough so that $\mu+\epsilon\in \Delta_K(M)$ is a weakly-regular value of $\Phi_K$.
\end{defi}

Let $n\geq 1$. The manifold $M$, equipped with the symplectic structure $n\Omega$,  is prequantized by the line bundle
$\Lcal^{\otimes n}$: the corresponding moment map is $n\Phi_K$. For any dominant weight $\mu \in  \wedge^*_{+}$, the symplectic reduction of 
$(M,n\Omega)$ relatively to the weight $n\mu$ is $(M_\mu,n\Omega_\mu)$. Like in definition \ref{def:RR-numbers}, we consider the following Riemann-Roch numbers
$$
\Qcal(M_\mu,n\Omega_\mu)=
\begin{cases}
0 \hspace{36mm} {\rm if }\quad  \mu\notin \Delta_K(M),\\
RR(M_{\mu+\epsilon},(\Lcal_{\mu,\epsilon})^{\otimes n})\hspace{4mm}  {\rm if }\quad \mu\in \Delta_K(M)\ {\rm and} \  \|\epsilon\|<<1.
\end{cases}
$$

Kawasaki-Riemann-Roch formula shows that $n\geq 1\mapsto \Qcal(M_\mu,n\Omega_\mu)$ is a quasi-polynomial map 
\cite{vergne96,Loizides19}. When $\mu$ is a weakly-regular value of $\Phi_K$, we denote by 
${\rm vol}(M_\mu):=\frac{1}{d_\mu}\int_{M_\mu}(\frac{\Omega_\mu}{2\pi})^{\frac{\dim M_\mu}{2}}$ the symplectic volume of the symplectic orbifold 
$(M_\mu,\Omega_\mu)$. Here $d_\mu$ is the generic value of the map $m\in \Phi_K^{-1}(\mu)\mapsto {\rm cardinal}(K_m/K^0_m)$.

The following proposition is a direct consequence of Kawasaki-Riemann-Roch formula (see \cite{Loizides19} or \S 1.3.4 in \cite{pep-duistermaat}).

\begin{prop}\label{prop:quasi-polynomial-non-zero}
Let $\mu\in\Delta_K(M)\cap \wedge^*_{+}$  be a weakly-regular value of $\Phi_K$. Then we have
$\Qcal(M_\mu,n\Omega_\mu)\ \sim\ {\rm vol}(M_\mu) \,n^{\frac{\dim M_\mu}{2}}$, when $n\to\infty$. 
In particular the map $n\geq 1\mapsto \Qcal(M_\mu,n\Omega_\mu)$ is non-zero.
\end{prop}

\medskip

\subsection{Quantization commutes with reduction}\label{sec:QR=0}

Let us explain the  ``{\em Quantization commutes with Reduction}'' Theorem proved in \cite{pep-jems}. 
 
We fix $\tilde{\lambda}\in \wGhol$. The coadjoint orbit $\wG\tilde{\lambda}$ is prequantized by the line bundle 
 $\wG\times_{K_{\tilde{\lambda}}}\Cbb_{\tilde{\lambda}}$, and the moment map $\Phi_G^{\tilde{\lambda}}:\wG\tilde{\lambda}\to\ggot^*$ 
 corresponding to the $G$-action on 
 $\wG\times_{K_{\tilde{\lambda}}}\Cbb_{\tilde{\lambda}}$ is equal to the restriction of the map $\pi_{\ggot,\tilde{\ggot}}$ to $\wG\tilde{\lambda}$.
 
The symplectic slice $Y_{\tilde{\lambda}}=(\Phi_G^{\tilde{\lambda}})^{-1}(\kgot^*)$ is prequantized by the line bundle 
$\Lcal_{\tilde{\lambda}}:=\wG\times_{K_{\tilde{\lambda}}}\Cbb_{\tilde{\lambda}}\vert_{Y_{\tilde{\lambda}}}$. The moment map 
$\Phi_K^{\tilde{\lambda}}:Y_{\tilde{\lambda}}\to\kgot^*$ corresponding to the $K$-action  is equal to the restriction of 
$\Phi_G^{\tilde{\lambda}}$ to $Y_{\tilde{\lambda}}$.

For any $\lambda\in \Ghol$, we consider the (possibly singular) symplectic reduced space
$$
\mathbb{X}_{\tilde{\lambda},\lambda}:=(\Phi_K^{\tilde{\lambda}})^{-1}(\lambda)/K_\lambda,
$$
equipped with the reduced symplectic form $\Omega_{\tilde{\lambda},\lambda}$, and the (possibly singular) line bundle
$$
\mathbb{L}_{\tilde{\lambda},\lambda}:=\left(\Lcal_{\tilde{\lambda}}\vert_{(\Phi_K^{\tilde{\lambda}})^{-1}(\lambda)}\otimes \Cbb_{-\lambda}\right)/K_\lambda.
$$

Thanks to Definition \ref{def:RR-numbers}, the geometric quantization 
$\Qcal(\mathbb{X}_{\tilde{\lambda},\lambda},\Omega_{\tilde{\lambda},\lambda})\in\Zbb$ of those compact symplectic spaces 
$(\mathbb{X}_{\tilde{\lambda},\lambda},\Omega_{\tilde{\lambda},\lambda})$ are well-defined even if they are singular. 
In particular $\Qcal(\mathbb{X}_{\tilde{\lambda},\lambda},\Omega_{\tilde{\lambda},\lambda})=0$ when 
$\mathbb{X}_{\tilde{\lambda},\lambda}:=\emptyset$.

The following Theorem is proved in \cite{pep-jems}.

\begin{theo}\label{theo:QR-jems}
Let $\tilde{\lambda}\in \wGhol$. We have an Hilbertian direct sum 
$$
V^{\wG}_{\tilde{\lambda}}\vert_{G}=\bigoplus_{\lambda\in \Ghol} \Qcal(\mathbb{X}_{\tilde{\lambda},\lambda},\Omega_{\tilde{\lambda},\lambda}) \ V^{G}_\lambda,
$$
It means that for any $\lambda\in \Ghol$, the multiplicity of the representation $V^{G}_\lambda$ in 
the restriction $V^{\wG}_{\tilde{\lambda}}\vert_{G}$ is equal to the geometric quantization
$\Qcal(\mathbb{X}_{\tilde{\lambda},\lambda},\Omega_{\tilde{\lambda},\lambda})$ of the (compact) reduced space $\mathbb{X}_{\tilde{\lambda},\lambda}$.
\end{theo}

\begin{rem}
Let $({\tilde{\lambda},\lambda})\in \wGhol\times\Ghol$. Theorem \ref{theo:QR-jems}. shows that 
$$
\left[V^{G}_{n\lambda} : V^{\wG}_{n\tilde{\lambda}}\right]= \Qcal(\mathbb{X}_{\tilde{\lambda},\lambda},n\Omega_{\tilde{\lambda},\lambda})
$$
for any $n\geq 1$.
\end{rem}

\section{Proofs of the main results}\label{sec:proof}

We come back to the setting of \S \ref{sec:cone-closed}:  $G/K$ is a complex submanifold of a Hermitian symmetric space $\wG/\wK$. 
It means that there exits a $\tilde{K}$-invariant element $z\in\kgot$ such that $\ad(z)$ defines complex structures on $\tilde{\pgot}$ and $\pgot$.
We consider the orthogonal decomposition $\tilde{\pgot}=\pgot\oplus\qgot$, and we denote by ${\rm Sym}(\qgot)$ the symmetric algebra of the complex $K$-module 
$(\qgot,{\rm ad}(z))$.

\subsection{Proof of Theorem A}\label{sec:proof-theo-A}
The set $\Pi_{\rm hol}(\tilde{G},G)$ is equal to $\bigcup_{\tilde{a}\in \wchol}\{\tilde{a}\}\times\Delta_G(\wG\tilde{a})$. We define
$$
\Pi_{\rm hol}(\tilde{G},G)^0:=\bigcup_{\tilde{a}\in \wchol} \{\tilde{a}\} \times \Delta_G(\wG\tilde{a})^0.
$$

We start with the following result.

\begin{lem}\label{lem:dense-Pi}
The set $\Pi_{\rm hol}(\tilde{G},G)^0\bigcap \tilde{\tgot}^*_\Qbb \times \tgot^*_\Qbb$ is dense in $\Pi_{\rm hol}(\tilde{G},G)$.
\end{lem}
{\em Proof :} 
Let $(\tilde{\xi},\xi)\in \Pi_{\rm hol}(\tilde{G},G)$: take $\tilde{g}\in\wG$ such that $\xi=\pi_{\ggot,\tilde{\ggot}}(\tilde{g}\tilde{\xi})$. We consider a sequence 
$\tilde{\xi}_n\in \wchol\cap\tilde{\tgot}^*_\Qbb$ converging to $\tilde{\xi}$. Then $\xi_n:=\pi_{\ggot,\tilde{\ggot}}(\tilde{g}\tilde{\xi}_n)$ is a sequence of $\Ccal^0_{G/K}$ converging to 
$\xi\in\chol$. Since the map ${\rm p}: \Ccal_{G/K}^0\to \chol$ is continuous (se Lemma \ref{lem:p-continu}), the sequence $\eta_n:={\rm p}(\xi_n)$ converges to 
${\rm p}(\xi)=\xi$. By definition, we have $\eta_n\in \Delta_G(\wG\,\tilde{\xi}_n)$ for any $n\in\Nbb$. Since $\tilde{\xi}_n$ are rational, each subset 
$\Delta_G(\wG\tilde{\xi}_n)^0\cap \tgot_\Qbb^*$ is dense in $\Delta_G(\wG\tilde{\xi}_n)$ (see Lemma \ref{lem:rational-dense-coadjoint}). Hence, $\forall n\in\Nbb$, there exists $\zeta_n\in \Delta_G(\wG\,\tilde{\xi}_n)^0\cap \tgot^*_\Qbb$ such that 
$\|\zeta_n-\eta_n\|\leq 2^{-n}$. Finally, we see that $(\tilde{\xi}_n,\zeta_n)$ is a sequence of rational elements of $\Pi_{\rm hol}(\tilde{G},G)^0$ converging to 
$(\xi,\tilde{\xi})$.
$\Box$

\medskip

The main purpose of this section is the proof of the following theorem.

\medskip

\begin{theo}\label{theo:rational-Pi}
The following inclusions hold
\begin{eqnarray*}
\Pi_{\rm hol}(\tilde{G},G)^0\, \bigcap\, \tilde{\tgot}^*_\Qbb \times \tgot^*_\Qbb
\quad \underset{(1)}{\mbox{\Large$\subset$}}\quad 
\Pi_{\rm hol}^\Qbb(\tilde{G},G)
\quad \underset{(2)}{\mbox{\Large$\subset$}}\quad 
\Pi_{\rm hol}(\tilde{G},G).
\end{eqnarray*}
\end{theo}

Lemma \ref{lem:dense-Pi} and Theorem \ref{theo:rational-Pi} gives the important Corollary.

\begin{coro}\label{coro: pi-hol-dense}
$\Pi_{\rm hol}^\Qbb(\tilde{G},G)$ is dense in $\Pi_{\rm hol}(\tilde{G},G)$.
\end{coro}

\medskip

{\em Proof of Theorem \ref{theo:rational-Pi} :} 
Let $(\tilde{\mu},\mu)\in \Pi_{\rm hol}^\Qbb(\tilde{G},G)$: there exists $N\geq 1$ such that 
$(N\tilde{\mu},N\mu)\in\Pi_{\rm hol}^\Zbb(\tilde{G},G)$. The multiplicity $[V^{G}_{N\mu}:V^{\wG}_{N\tilde{\mu}}]$ is non-zero, and thanks to Theorem 
\ref{theo:QR-jems}, it implies that the reduced space $\mathbb{X}_{N\tilde{\mu},N\mu}$ is non-empty. In other words 
$(N\tilde{\mu},N\mu)\in \Pi_{\rm hol}(\tilde{G},G)$. The inclusion $(2)$ is proven.

Let $(\tilde{\mu},\mu)\in \Pi_{\rm hol}(\tilde{G},G)^0\bigcap \tgot^*_\Qbb\times\tilde{\tgot}^*_\Qbb$. There exists $N_o\geq 1$ such that 
$\lambda:=N_o\mu \in \Ghol$,  $\tilde{\lambda}:=N_o\tilde{\mu}\in\wGhol$ and $\lambda\in \Delta_G(\wG\tilde{\lambda})^0$: the element $\lambda$ is a weakly-regular value of the moment map $\wG\tilde{\lambda}\to\ggot^*$. Theorem \ref{theo:QR-jems} tell us that, for any $n\geq 1$,  the multiplicity 
$[V^{G}_{n\lambda}:V^{\wG}_{n\tilde{\lambda}}]$ is equal to Riemann-Roch number 
$\Qcal(\mathbb{X}_{\tilde{\lambda},\lambda}, n\Omega_{\tilde{\lambda},\lambda})$. Since the map 
$n\mapsto \Qcal(\mathbb{X}_{\tilde{\lambda},\lambda}, n\Omega_{\tilde{\lambda},\lambda})$ is non-zero 
(see Proposition \ref{prop:quasi-polynomial-non-zero}), we can conclude that there exists $n_o\geq 1$ such that 
$[V^{G}_{n_o\lambda}:V^{\wG}_{n_o\tilde{\lambda}}]\neq 0$. In other words,  we obtain $n_oN_o(\tilde{\mu},\mu)\in \Pi_{\rm hol}^\Zbb(\tilde{G},G)$ 
and so $(\tilde{\mu},\mu)\in \Pi_{\rm hol}^\Qbb(\tilde{G},G)$. The inclusion $(1)$ is settled.
$\Box$

\medskip

\medskip

Now we can terminate the proof of Theorem A.

\medskip

Thanks to Proposition \ref{prop:monoid}, we know that $\Pi_{\rm hol}^\Qbb(\tilde{G},G)$ is a $\Qbb$-convex cone.  Since $\Pi_{\rm hol}(\tilde{G},G)$ is a closed subset of $\wchol\times\chol$ (see Proposition \ref{prop:pi-z-compact}), 
we can conclude, by a density argument,  that $\Pi_{\rm hol}(\tilde{G},G)$ is a closed convex cone of $\wchol\times\chol$.

\subsection{The affine variety $\wK_\Cbb\times\qgot$}\label{sec:affine}

Let $\tilde{\kappa}$ be the Killing form on the Lie algebra $\tilde{\ggot}$. We consider the $\wK$-invariant symplectic structures $\Omega_{\tilde{\pgot}}$ on $\tilde{\pgot}$, defined by the relation 
$$
\Omega_{\tilde{\pgot}}(\tilde{Y},\tilde{Y}')=\tilde{\kappa}(z,[\tilde{Y},\tilde{Y}']),\quad \forall \tilde{Y},\tilde{Y}'\in\tilde{\pgot}.
$$
We notice that the complex structure ${\rm ad}(z)$ is adapted to $\Omega_{\tilde{\pgot}}$ : $\Omega_{\tilde{\pgot}}(\tilde{Y},{\rm ad}(z)\tilde{Y})>0$ if $\tilde{Y}\neq 0$.

We denote by $\Omega_{\qgot}$ the restriction of 
$\Omega_{\tilde{\pgot}}$ on the symplectic subspace $\qgot$. The moment map $\Phi_{\qgot}$ associated to the $K$-action on 
$(\qgot, \Omega_{\qgot})$ is defined by the relations $\langle\Phi_{\qgot}(Y),X\rangle= \frac{-1}{2}\tilde{\kappa}([X,Y],[z,Y])$,  
$\forall (X,Y) \in\pgot\times \qgot$.  In particular, $\langle\Phi_{\qgot}(Y),z\rangle = \frac{-1}{2}\|Y\|^2$, $\forall Y \in\qgot$, 
so the map $\langle\Phi_{\qgot},z\rangle$ is proper.

The complex reductive group $\tilde{K}_\Cbb$ is equipped with the following action of 
$\tilde{K}\times K$ : $(\tilde{k},k)\cdot a= \tilde{k} a k^{-1}$.   It has a canonical structure of a smooth affine variety. 
There is a diffeomorphism of the cotangent bundle $\T^*\tilde{K}$ with $\tilde{K}_\Cbb$ defined as follows. We identify $\T^*\tilde{K}$ with $\tilde{K}\times \tilde{\kgot}^*$ 
by means of left-translation and then with $\tilde{K}\times \tilde{\kgot}$ by means of an invariant inner product on  $\tilde{\kgot}$. The map $\varphi:\tilde{K}\times \tilde{\kgot}\to \tilde{K}_\Cbb$ given by $\varphi(a,X)=a e^{iX}$ is a diffeomorphism. If we use $\varphi$ to transport the complex structure of $\tilde{K}_\Cbb$ to $\T^*\tilde{K}$, then the resulting complex structure on $\T^*\tilde{K}$ is compatible with the symplectic structure on $\T^*\tilde{K}$, so that $\T^*\tilde{K}$ 
becomes a K\"{a}hler Hamiltonian $\tilde{K}\times K$-manifold (see \cite{Hall97}, \S 3). The moment map relative to the $\tilde{K}\times K$-action is the proper 
map $\Phi_{\tilde{K}} \oplus \Phi_K : \T^*\tilde{K} \to \tilde{\kgot}^*\oplus \kgot^*$ defined by
$\Phi_{\tilde{K}}(\tilde{a},\tilde{\eta})=-\tilde{a}\tilde{\eta}$ and $\Phi_K(\tilde{a},\tilde{\eta})=\pi_{\kgot,\wk}(\tilde{\eta})$. 
Here $\pi_{\kgot,\wk}: \tilde{\kgot}^* \to \kgot^*$ is the projection dual to the inclusion $\kgot \croc \tilde{\kgot}$ of Lie algebras.

Finally we consider the K\"{a}hler Hamiltonian $\tilde{K}\times K$-manifold  
$\T^*\tilde{K}\times \qgot$, where $\qgot$ is equipped with the symplectic structure $\Omega_\qgot$. 
Let us denote by $\Phi: \T^*\tilde{K}\times \qgot\to \tilde{\kgot}^*\oplus \kgot^*$ 
the moment map relative to the $\tilde{K}\times K$-action~:
\begin{equation}\label{eq:momentcotangent}
\Phi(\tilde{a},\tilde{\eta},Y)=\left(-\tilde{a}\tilde{\eta},\pi_{\kgot,\wk}(\tilde{\eta})+\Phi_\qgot(Y)\right).
\end{equation}

Since $\Phi$ is proper map, the Convexity Theorem tell us that 
$$
\Delta(\T^*\tilde{K}\times \qgot):={\rm Image}(\Phi)\bigcap \, \tilde{\tgot}^*_{\geq 0}\times \tgot^*_{\geq 0}
$$
is a closed convex locally polyhedral set.

We consider now the action of $\wK\times K$ on the affine variety $\tilde{K}_\Cbb\times \qgot$. The set of highest weights of 
$\tilde{K}_\Cbb\times \overline{\qgot}$ is the semigroup 
$$
\Delta^{\Zbb}(\tilde{K}_\Cbb\times\qgot)\subset \tilde{\wedge}^*_+\times\wedge^*_+
$$ 
consisting of all dominant weights $(\tilde{\lambda},\lambda)$ such that the irreducible $\wK\times K$-representation 
$V_{\tilde{\lambda}}^{\wK}\otimes V_{\lambda}^{K}$ occurs in the coordinate ring $\Cbb[\tilde{K}_\Cbb\times \qgot]$. 
A direct application of the Peter-Weyl Theorem gives the following characterization : 
\begin{equation}\label{eq:delta-Z-affine}
(\tilde{\lambda},\lambda)\in\Delta^{\Zbb}(\tilde{K}_\Cbb\times \qgot)\Longleftrightarrow
\left[ V^{\wK}_{\tilde{\lambda}}\vert_K\otimes V^{K}_\lambda\otimes {\rm Sym}(\qgot) \right]^K\neq 0.
\end{equation}

We denote by $\Delta^{\Qbb}(\tilde{K}_\Cbb\times \qgot)$ the $\Qbb$-convex cone generated by the semigroup \break $\Delta^{\Zbb}(\tilde{K}_\Cbb\times \qgot)$: 
$(\tilde{\xi},\xi)\in \Delta^{\Qbb}(\tilde{K}_\Cbb\times \qgot)$ if and only if $\exists N\geq 1$, $N(\tilde{\xi},\xi)\in \Delta^{\Zbb}(\tilde{K}_\Cbb\times \qgot)$.

\medskip

The following important fact is classical (see Theorem 4.9 in \cite{sjamaar98}).

\begin{prop}\label{prop:delta-affine-dense}
The Kirwan polyhedron $\Delta(\T^*\tilde{K}\times \qgot)$ is equal to the closure of the $\Qbb$-convex cone $\Delta^{\Qbb}(\tilde{K}_\Cbb\times \qgot)$.
\end{prop}

\subsection{Proof of Theorem B}
Consider the semigroup $\Pi^\Zbb(\tilde{G},G)$ of $\tilde{\wedge}^*_+ \times \wedge^*_+$ (see Definition \ref{defi:horn-cone-general}) and the $\Qbb$-convex cone  
$\Pi^\Qbb(\tilde{G},G):=\{(\tilde{\xi},\xi)\in \tilde{\tgot}^*_{\geq 0}\times \tgot^*_{\geq 0} \ ; \ \exists N\geq 1, N(\tilde{\xi},\xi)\in \Pi^\Zbb(\tilde{G},G)\}$.

The Jakobsen-Vergne theorem says that $\Pi_{\rm hol}^\Zbb(\tilde{G},G)=\Pi^\Zbb(\tilde{G},G)\,\bigcap \, \wGhol\times\Ghol$. Hence, the convex cone $\Pi_{\rm hol}^\Qbb(\tilde{G},G)$ is equal to 
$\Pi^\Qbb(\tilde{G},G)\cap \wchol\times\chol$. Thanks to (\ref{eq:delta-Z-affine}), we know that 
$(\tilde{\xi},\xi)\in \Pi^\Qbb(\tilde{G},G)$ if and only if $(\tilde{\xi},\xi^*)\in \Delta^{\Qbb}(\tilde{K}_\Cbb\times \qgot)$. 
The density results obtained in Proposition \ref{prop:delta-affine-dense} and Corollary \ref{coro: pi-hol-dense} gives finally Theorem B.

\begin{theo}\label{theo:theorem-B}
\begin{enumerate}
\item We have $\Pi_{\rm hol}(\tilde{G},G)=\Pi(\tilde{G},G)\,\bigcap \,\wchol\times\chol $.
\item An element $(\tilde{\xi},\xi)$ belongs to  $\Pi_{\rm hol}(\wG,G)$ if and only if 
$$
(\tilde{\xi},\xi)\in\wchol\times\chol\quad {\rm and}\quad (\tilde{\xi},\xi^*)\in \Delta(\T^*\wK\times \qgot).
$$
\end{enumerate}
\end{theo}

\medskip

\subsection{Proof of Theorem C}

We denote by $\bar{\qgot}$ the $K$-vector space $\qgot$ equipped with the opposite symplectic form $-\Omega_{\qgot}$ and opposite complex structure 
$-{\rm ad}(z)$. The moment map relative to the $K$-action on $\bar{\qgot}$ is denoted by $\Phi_{\bar{\qgot}}=-\Phi_{\qgot}$.

\begin{lem} Any element $(\tilde{\xi},\xi)\in\tilde{\tgot}^*_{\geq 0}\times\tgot^*_{\geq 0}$ satisfies the equivalence : 
$$
(\tilde{\xi},\xi^*)\in \Delta(\T^*\wK\times \qgot) \Longleftrightarrow \xi\in \Delta_K(\wK\tilde{\xi}\times \overline{\qgot}).
$$
\end{lem}
{\em Proof :} Thanks to (\ref{eq:momentcotangent}), we see immediatly that $\exists (\tilde{a},\tilde{\eta},Y)\in\T ^*\wK\times \qgot$ such that 
$(\tilde{\xi},\xi^*)=\Phi(\tilde{a},\tilde{\eta},Y)$ if and only if $\exists (\tilde{b},Z)\in \wK\times \qgot$ such that 
$\xi=\pi_{\kgot,\wk}(\tilde{b}\tilde{\xi})+\Phi_{\bar{\qgot}}(Z)$. $\Box$

\medskip

At this stage we know that  $\Delta_G(\wG\tilde{\mu})=\Delta_K(\wK\tilde{\mu}\times \overline{\qgot})\cap \chol$. Hence, Theorem C will follows from the next result.

\begin{prop}\label{prop:delta-chol}
For any $\tilde{\mu}\in \wchol$, the Kirwan polyhedron $\Delta_K(\wK\tilde{\mu}\times \overline{\qgot})$ is contained in $\chol$.
\end{prop}

{\em Proof :}  By definition $\chol=\Ccal_{G/K}^0\cap\tgot^*_{\geq 0}$, so we have to prove that 
$\pi_{\kgot,\tilde{\kgot}}(\wK\tilde{\mu})+ {\rm Image}( \Phi_{\bar{\qgot}})$ is contained in $\Ccal_{G/K}^0$. 
By definition $\wK\tilde{\mu}\subset \Ccal_{\wG/\wK}^0$, and then $\pi_{\kgot,\tilde{\kgot}}(\wK\tilde{\mu})\subset \Ccal_{G/K}^0$. Since 
$\Ccal_{G/K}^0 +\Ccal_{G/K}\subset \Ccal_{G/K}^0$, it is sufficient to check that ${\rm Image}( \Phi_{\bar{\qgot}})\subset \Ccal_{G/K}$. 
Let $\Phi_{\tilde{\pgot}}$ be the moment map relative to the action of $\wK$ on $(\tilde{\pgot},\Omega_{\tilde{\pgot}})$. As 
${\rm Image}( \Phi_{\bar{\qgot}})\subset \pi_{\kgot,\tilde{\kgot}}\left({\rm Image}( -\Phi_{\tilde{\pgot}})\right)$, the following lemma will terminate the proof of 
Proposition \ref{prop:delta-chol}. $\Box$

\begin{lem}
The image of the moment map $-\Phi_{\tilde{\pgot}}$ is contained in $\Ccal_{\wG/\wK}$.
\end{lem}
{\em Proof :}  
Let $z^*\in\tilde{\tgot}^*$ such that $\langle z^*,\tilde{X}\rangle= -\tilde{\kappa}(z,\tilde{X})$, $\forall\tilde{X}\in\tilde{\ggot}$. Consider the coadjoint orbit 
$\tilde{\Ocal}=\wG z^*$ equipped with its canonical symplectic structure $\Omega_{\tilde{\Ocal}}$: the symplectic vector space $\T_{z^*}\tilde{\Ocal}$ is canonically 
isomorphic to $(\tilde{\pgot},-\Omega_{\tilde{\pgot}})$. In \cite{McDuff}, McDuff proved that $(\tilde{\Ocal},\Omega_{\tilde{\Ocal}})$ is diffeomorphic, as a 
$\wK$-symplectic manifold, to the symplectic vector space $(\tilde{\pgot},-\Omega_{\tilde{\pgot}})$ (see \cite{Deltour-these,Deltour-JDG} for a generalization of this fact). 
McDuff's results shows in particular that ${\rm Image}(-\Phi_{\tilde{\pgot}})=\pi_{\tilde{\ggot},\tilde{\kgot}}(\tilde{\Ocal})$. Our proof is completed if we check that 
$\pi_{\tilde{\ggot},\tilde{\kgot}}(\tilde{\Ocal})\subset\Ccal_{\wG/\wK}$: in other words, if $\langle\pi_{\tilde{\ggot},\tilde{\kgot}}(\tilde{g}_0\, z^*),\tilde{g}_1z\rangle\geq 0$, 
$\forall \tilde{g}_0,\tilde{g}_1\in \wG$. But
\begin{eqnarray*}
2\langle\pi_{\tilde{\ggot},\tilde{\kgot}}(\tilde{g}_0\, z^*),\tilde{g}_1\,z\rangle &=& \langle\tilde{g}_0\, z^*,\tilde{g}_1z+\Theta(\tilde{g}_1)z\rangle\\
&=& -\tilde{\kappa}(z,\tilde{g}_0^{-1}\tilde{g}_1\, z)-\tilde{\kappa}(z,\tilde{g}_0^{-1}\Theta(\tilde{g}_1)z).
\end{eqnarray*}
With (\ref{eq:Gz}) in hand, it is not difficult to see that $-\tilde{\kappa}(z,\tilde{g}\, z)\geq 0$ for every $\tilde{g}\in\wG$. We thus verified that 
$\pi_{\tilde{\ggot},\tilde{\kgot}}(\tilde{\Ocal})\subset\Ccal_{\wG/\wK}$. $\Box$

\section{Facets of the cone $\Pi_{\rm hol}(\tilde{G},G)$}\label{sec:cone:faces}

We come back to the framework of \S \ref{sec:affine}. We consider the K\"{a}hler Hamiltonian $\tilde{K}\times K$-manifold  
$\T^*\tilde{K}\times \qgot$. The moment map, $\Phi: \T^*\tilde{K}\times \qgot\to \tilde{\kgot}^*\oplus \kgot^*$, relative to the $\tilde{K}\times K$-action is defined by (\ref{eq:momentcotangent}).

In this section, we adapt to our case the result of \S 6 of \cite{pep-ressayre-pair} concerning the parametrization of the facets of Kirwan polyhedrons in 
terms of Ressayre's data.

\subsection{Admissible elements}
We choose maximal torus $\wT\subset\wK$ and $T\subset K$ such that $T\subset \wT$. Let $\Rgot_o$ and $\Rgot$ be respectively the set of roots for the action of $T$ on 
$(\wg/\ggot)\otimes\Cbb$ and $\ggot\otimes\Cbb$. Let $\tilde{\Rgot}$ be the set of roots for the action of $\wT$ on $\wg\otimes\Cbb$. 
Let $\Rgot^+\subset \Rgot$ and  $\tilde{\Rgot}^+\subset \tilde{\Rgot}$  be the systems of positive roots defined in (\ref{eq:R-plus}).
Let $W,\tilde{W}$ be the Weyl groups of $(T,K)$ and $(\wT,\wK)$. Let $w_o\in W$ be the longest element.

We start by introducing the notion of admissible elements. The group $\hom(U(1),T)$ admits a natural identification with the lattice 
$\wedge:=\frac{1}{2\pi}\ker(\exp:\tgot\to T)$. A vector $\gamma\in\tgot$ is called rational if it belongs to the $\Qbb$-vector space 
$\tgot_\Qbb$ generated by $\wedge$.

We consider the $\wK\times K$-action on $N:=\T^* \wK \times \qgot$. We associate to any subset $\Xcal\subset N$, the integer 
$\dim_{\wK\times K}(\Xcal)$ (see (\ref{eq:dim-K})). 

\begin{defi}
A non-zero element $(\tilde{\gamma},\gamma)\in\tilde{\tgot}\times \tgot$ 
is called {\em admissible} if the elements $\tilde{\gamma}$ and $\gamma$ are rational and if $\dim_{\wK\times K}(N^{(\tilde{\gamma},\gamma)})-\dim_{\wK\times K}(N)\in\{0,1\}$.
\end{defi}

If $\gamma\in\tgot$, we denote by $\Rgot_o\cap\gamma^\perp$ the subsets of weight vanishing against $\gamma$. We start with the 
following lemma whose proof is left to the reader (see \S 6.1.1 of \cite{pep-ressayre-pair}).

\begin{lem}
\begin{enumerate}
\item $N^{(\tilde{\gamma},\gamma)}\neq \emptyset$ if and only if $\tilde{\gamma}\in \tilde{W}\gamma$.
\item $\dim_{\wK\times K}(N)=\dim_T(\tilde{\ggot}/\ggot)=\dim(\tgot)-\dim({\rm Vect}(\Rgot_o))$.
\item For any $\tilde{w}\in\tilde{W}$, $\dim_{\wK\times K}(N^{(\tilde{w}\gamma,\gamma)})=
\dim_T(\tilde{\ggot}^\gamma/\ggot^\gamma)=\dim(\tgot)-\dim({\rm Vect}(\Rgot_o\cap\gamma^\perp))$.
\end{enumerate}
\end{lem}

The next result is a direct consequence of the previous lemma.

\begin{lem}
The admissible elements relative to the $\wK\times K$-action on $\T^* \wK \times \qgot$ are of the form $(\tilde{w}\gamma,\gamma)$ where 
$\tilde{w}\in \tilde{W}$ and $\gamma$ is 
a non-zero rational element satisfying ${\rm Vect}(\Rgot_o)\cap \gamma^\perp={\rm Vect}(\Rgot_o\cap\gamma^\perp)$.
\end{lem}

\subsection{Ressayre's data}

\begin{defi}
\begin{enumerate}
\item Consider the linear action $\rho: G\to {\rm GL}_\Cbb(V)$ of a compact Lie group on a complex vector space $V$. For any $(\eta,a)\in\ggot\times\Rbb$, we define the vector subspace
$V^{\eta=a}=\{v\in V, d\rho(\eta)v=i av\}$.  Thus, for any $\eta\in\ggot$, we have the decomposition
$V=V^{\eta>0}\oplus V^{\eta=0}\oplus V^{\eta<0}$ where $V^{\eta>0}=\sum_{a>0}V^{\eta=a}$, and $V^{\eta<0}=\sum_{a<0}V^{\eta=a}$.
\item The real number $\tr_{\eta}(V^{\eta>0})$ is defined as the sum $\sum_{a>0}a\,\dim(V^{\eta=a})$.
\end{enumerate}
\end{defi}

We consider an admissible element $(\tilde{w}\gamma,\gamma)$. The submanifold of $N\simeq \wK_\Cbb \times \qgot$ fixed  of by $(\tilde{w}\gamma,\gamma)$ 
is $N^{(\tilde{w}\gamma,\gamma)}= \tilde{w} \wK_\Cbb^\gamma \times \qgot^{\gamma}$. There is a canonical isomorphism between the manifold 
$N^{(\tilde{w}\gamma,\gamma)}$ equipped with the action of $\tilde{w}\wK^\gamma \tilde{w}^{-1}\times K^\gamma$ 
with the manifold $\wK_\Cbb^\gamma \times \qgot^{\gamma}$ equipped with the action of $\wK^\gamma \times K^\gamma$.
The tangent bundle $(\T N\vert_{N^{(\tilde{w}\gamma,\gamma)}})^{(\tilde{w}\gamma,\gamma)>0}$ is isomorphic to 
$N^{\gamma_w}\times \tilde{\kgot}_\Cbb^{\gamma>0}\times \qgot^{\gamma>0}$. 

The choice of positive roots $\Rgot^+$ (resp. $\tilde{\Rgot}^+$) induces a decomposition $\kgot_\Cbb=\ngot\oplus\tgot_\Cbb\oplus \overline{\ngot}$ (resp. $\wk_\Cbb=\tilde{\ngot}\oplus\tilde{\tgot}_\Cbb\oplus \overline{\tilde{\ngot}}$), where 
$\ngot=\sum_{\alpha\in\Rgot^+}(\kgot\otimes\Cbb)_\alpha$ (resp. $\tilde{\ngot}=\sum_{\tilde{\alpha}\in\tilde{\Rgot}^+}(\wk\otimes\Cbb)_{\tilde{\alpha}}$).
 We consider the map 
$$
\rho^{\tilde{w},\gamma}: \wK_\Cbb^\gamma \times \qgot^\gamma \longrightarrow 
\hom \left(\tilde{\ngot}^{\tilde{w}\gamma>0}\times \ngot^{\gamma>0}\,,\, \tilde{\kgot}_\Cbb^{\gamma>0}\times\qgot^{\gamma>0}\right)
$$
defined by the relation
$$
\rho^{\tilde{w},\gamma}(\tilde{x},v): (\tilde{X},X)\longmapsto ((\tilde{w} \tilde{x})^{-1}\tilde{X}-X\, ;\, X\cdot v),
$$
for any $(\tilde{x},v)\in \wK_\Cbb^\gamma \times \qgot^\gamma$.
\medskip

\begin{defi}
$(\gamma,\tilde{w})\in\tgot\times \tilde{W}$ is a Ressayre's data if 
\begin{enumerate}
\item $(\tilde{w}\gamma,\gamma)$ is admissible,
\item $\exists (\tilde{x},v)$ such that $\rho^{\tilde{w},\gamma}(\tilde{x},v)$ is bijective.
\end{enumerate}
\end{defi}

\begin{rem}
In \cite{pep-ressayre-pair}, the Ressayre's data were called {\em regular infinitesimal $B$-Ressayre's pairs}.
\end{rem}

Since the linear map $\rho^{\tilde{w},\gamma}(\tilde{x},v)$ commutes with the $\gamma$-actions, we obtain the following necessary conditions.

\begin{lem}\label{lem:relations-A-B}
If $(\gamma,\tilde{w})\in\tgot\times \tilde{W}$ is a Ressayre's data, then
\begin{itemize}
\item Relation (A) : $\dim(\tilde{\ngot}^{\tilde{w}\gamma>0})+\dim(\ngot^{\gamma>0})=\dim(\tilde{\kgot}_\Cbb^{\gamma>0})+\dim(\qgot^{\gamma>0})$.
\item Relation (B) : $\tr_{\tilde{w}\gamma}(\tilde{\ngot}^{\tilde{w}\gamma>0})+\tr_{\gamma}(\ngot^{\gamma>0})=
\tr_{\gamma}(\tilde{\kgot}_\Cbb^{\gamma>0})+\tr_{\gamma}(\qgot^{\gamma>0})$.
\end{itemize}
\end{lem}

\begin{lem}
Relation (B) is equivalent to 
\begin{equation}\label{eq:trace-condition-w-tilde}
\sum_{\stackrel{\alpha\in\Rgot^+}{\langle\alpha,\gamma\rangle> 0}}\langle\alpha,\gamma\rangle=
\sum_{\stackrel{\tilde{\alpha}\in\tilde{\Rgot}^+}{\langle\tilde{\alpha},\tilde{w}_0\tilde{w}\gamma\rangle> 0}}\langle\tilde{\alpha},\tilde{w}_0\tilde{w}\gamma\rangle.
\end{equation}
\end{lem}

{\em Proof :} First one sees that  
$\tr_{\gamma}(\qgot^{\gamma>0})=\tr_{\gamma}(\tilde{\pgot}^{\gamma>0})-\tr_{\gamma}(\pgot^{\gamma>0})=
\sum_{\stackrel{\tilde{\alpha}\in \tilde{\Rgot}^+_n}{\langle\tilde{\alpha},\gamma\rangle>0}}\langle\tilde{\alpha},\gamma\rangle
-\sum_{\stackrel{\alpha\in \Rgot^+_n}{\langle\alpha,\gamma\rangle>0}}\langle\alpha,\gamma\rangle$, and 
$\tr_{\gamma}(\wk_\Cbb^{\gamma>0})=\tr_{\tilde{w}\gamma}(\wk_\Cbb^{\tilde{w}\gamma>0})=\tr_{\tilde{w}\gamma}(\tilde{\ngot}^{\tilde{w}\gamma>0})+
\sum_{\stackrel{\tilde{\alpha}\in \tilde{\Rgot}^+_c}{\langle\tilde{\alpha},\tilde{w}_0\tilde{w}\gamma\rangle>0}}\langle\tilde{\alpha},\tilde{w}_0\tilde{w}\gamma\rangle$.
Relation $(B)$ is equivalent to 
\begin{equation}\label{eq:equivalent-relation-B}
\tr_{\gamma}(\ngot^{\gamma>0})+\sum_{\stackrel{\alpha\in \Rgot^+_n}{\langle\alpha,\gamma\rangle>0}}\langle\alpha,\gamma\rangle=
\sum_{\stackrel{\tilde{\alpha}\in \tilde{\Rgot}^+_n}{\langle\tilde{\alpha},\gamma\rangle>0}}\langle\tilde{\alpha},\gamma\rangle
+\sum_{\stackrel{\tilde{\alpha}\in \tilde{\Rgot}^+_c}{\langle\tilde{\alpha},\tilde{w}_0\tilde{w}\gamma\rangle>0}}\langle\tilde{\alpha},\tilde{w}_0\tilde{w}\gamma\rangle.
\end{equation}
Since $\tilde{\Rgot}^+_n$ is invariant under the action of the Weyl group $\tilde{W}$, the right hand side of (\ref{eq:equivalent-relation-B}) is equal to
$\sum_{\stackrel{\tilde{\alpha}\in\tilde{\Rgot}^+}{\langle\tilde{\alpha},\tilde{w}_0\tilde{w}\gamma\rangle>0}}\langle\tilde{\alpha},\tilde{w}_0\tilde{w}\gamma\rangle$.
Since the left hand side of (\ref{eq:equivalent-relation-B}) is equal to 
$\sum_{\stackrel{\alpha\in \Rgot^+}{\langle\alpha,\gamma\rangle>0}}\langle\alpha,\gamma\rangle$, the proof of the Lemma is complete. $\Box$

\subsection{Cohomological characterization of Ressayre's data}

Let $\gamma\in\tgot$ be a non-zero rational element. We denote by $B\subset K_\Cbb$ and by $\tilde{B}\subset \wK_\Cbb$ the Borel subgroups with 
Lie algebra $\bgot=\tgot_\Cbb\oplus \ngot$ and $\tilde{\bgot}=\tilde{\tgot}_\Cbb\oplus \tilde{\ngot}$. Consider the parabolic subgroup $P_\gamma\subset K_\Cbb$ defined by
\begin{equation}\label{eq:P-gamma}
P_\gamma=\{g\in K_\Cbb, \lim_{t\to\infty}\exp(-it\gamma)g\exp(it\gamma)\ {\rm exists}\}.
\end{equation}
Similarly, one defines a parabolic subgroup $\tilde{P}_\gamma\subset \wK_\Cbb$.

We work with the projective varieties $\Fcal_\gamma:=K_\Cbb/P_\gamma$, $\tilde{\Fcal}_\gamma:=\tilde{K}_\Cbb/\tilde{P}_\gamma$, 
and the canonical embedding $\iota: \Fcal_\gamma\to  \tilde{\Fcal}_\gamma$. We associate to any $\tilde{w}\in\tilde{W}$, the Schubert cell
$$
\tilde{\Xgot}^o_{\tilde{w},\gamma}:= \tilde{B}[\tilde{w}]\subset \tilde{\Fcal}_\gamma.
$$
and the Schubert variety $\tilde{\Xgot}_{\tilde{w},\gamma}:=\overline{\tilde{\Xgot}^o_{\tilde{w},\gamma}}$. If $\tilde{W}^\gamma$ denotes the subgroup of $\tilde{W}$ that fixes 
$\gamma$, we see that  the Schubert cell $\tilde{\Xgot}^o_{\tilde{w},\gamma}$ and the Schubert variety $\tilde{\Xgot}_{\tilde{w},\gamma}$ depends only of the class of $\tilde{w}$ in 
$\tilde{W}/\tilde{W}^\gamma$.

On the variety $\Fcal_\gamma$, we consider the Schubert cell $\Xgot^o_{\gamma}:= B[e]$ and the Schubert variety $\Xgot_{\gamma}:=\overline{\Xgot^o_{\gamma}}$. 

We consider the cohomology\footnote{Here, we use singular cohomology with integer coefficients.} ring $H^*(\tilde{\Fcal}_\gamma,\Zbb)$ of $\tilde{\Fcal}_\gamma$. 
If Y is an irreducible closed subvariety of $\tilde{\Fcal}_\gamma$, we denote by $[Y]\in H^{2n_Y}(\tilde{\Fcal}_\gamma,\Zbb)$ its cycle class in cohomology : here $n_Y={\rm codim}_\Cbb(Y)$. Let $\iota^*:H^*(\tilde{\Fcal}_\gamma,\Zbb)\to H^*(\Fcal_\gamma,\Zbb)$ be the pull-back map in cohomology. 
Recall that the cohomology class $[pt]$  associated to a singleton $Y=\{pt\}\subset \Fcal_\gamma$   is a basis of 
$H^{\mathrm{max}}(\Fcal_\gamma,\Zbb)$.

\medskip

To an oriented real vector bundle $\Ecal\to N$ of rank $r$, we can associate its Euler class ${\rm Eul}(\Ecal)\in H^{2r}(N,\Zbb)$. 
When $\Vcal\to N$ is a complex vector bundle, then ${\rm Eul}(\Vcal_\Rbb)$ corresponds to the top Chern class $c_{p}(\Vcal)$, where $p$ is the complex rank of $\Vcal$, and 
$\Vcal_\Rbb$ means $\Vcal$ viewed as a real vector bundle oriented by its complex structure (see \cite{Bott-Tu}, \S 21).

The isomorphism $\qgot^{\gamma>0}\simeq \qgot/ \qgot^{\gamma\leq 0}$ shows that $\qgot^{\gamma>0}$ can be viewed as a $P_\gamma$-module. Let 
$[\qgot^{\gamma>0}]= K_\Cbb\times_{P_\gamma}\qgot^{\gamma>0}$ be the corresponding complex vector bundle on $\Fcal_\gamma$.
We denote simply by ${\rm Eul}(\qgot^{\gamma>0})$ the Euler class ${\rm Eul}([\qgot^{\gamma>0}]_\Rbb)\in H^*(\Fcal_\gamma,\Zbb)$.

The following characterization of Ressayre's data was obtained in \cite{pep-ressayre-pair}, \S 6. Recall that $\Rgot_o$ denotes the set of weights relative to the $T$-action on 
$(\wg/\ggot)\otimes\Cbb$.

\begin{prop}
An element $(\gamma,\tilde{w})\in\tgot\times \tilde{W}$ is a Ressayre's data if and only if the following conditions hold :
\begin{itemize}
\item $\gamma$ is non-zero and rational.
\item ${\rm Vect}(\Rgot_o\cap \gamma^\perp)={\rm Vect}(\Rgot_o)\cap \gamma^\perp$.
\item $[\Xgot_{\gamma}]\cdot \iota^*([\tilde{\Xgot}_{\tilde{w},\gamma}])\cdot {\rm Eul}(\qgot^{\gamma>0})= k[pt],\ k\geq 1$ in $H^*(\Fcal_\gamma,\Zbb)$.
\item $\sum_{\stackrel{\alpha\in\Rgot^+}{\langle\alpha,\gamma\rangle> 0}}\langle\alpha,\gamma\rangle=
\sum_{\stackrel{\tilde{\alpha}\in\tilde{\Rgot}^+}{\langle\tilde{\alpha},\tilde{w}_0\tilde{w}\gamma\rangle> 0}}\langle\tilde{\alpha},\tilde{w}_0\tilde{w}\gamma\rangle$.
\end{itemize}
\end{prop}

\subsection{Parametrization of the facets}\label{sec:facets}

We can finally describe the Kirwan polyhedron $\Delta(\T^*\tilde{K}\times \qgot)$ (see \cite{pep-ressayre-pair}, \S 6).

\begin{theo}\label{theo:delta} 
An element $(\tilde{\xi},\xi)\in\tilde{\tgot}_{\geq 0}^*\times \tgot_{\geq 0}^*$ belongs to $\Delta(\T^*\tilde{K}\times  \qgot)$ if and only if
$$
\langle \tilde{\xi},\tilde{w}\gamma\rangle+\langle \xi,\gamma\rangle\geq 0
$$
for any Ressayre's data $(\gamma,\tilde{w})\in\tgot\times \tilde{W}$.
\end{theo}

Theorem \ref{theo:delta} and Theorem \ref{theo:theorem-B} permit us to give the following description of the convex cone 
$\Pi_{\rm hol}(\tilde{G},G)$.

\begin{theo}
An element $(\tilde{\xi},\xi)$ belongs to  $\Pi_{\rm hol}(\wG,G)$ if and only if $(\tilde{\xi},\xi)\in\wchol\times\chol$ and
$$
\langle \tilde{\xi},\tilde{w}\gamma\rangle\geq \langle \xi,w_0\gamma\rangle
$$
for any $(\gamma,\tilde{w})\in\tgot\times \tilde{W}$ satisfying the following conditions:
\begin{itemize}
\item $\gamma$ is non-zero and rational.
\item ${\rm Vect}(\Rgot_o\cap \gamma^\perp)={\rm Vect}(\Rgot_o)\cap \gamma^\perp$.
\item $[\Xgot_{\gamma}]\cdot \iota^*([\tilde{\Xgot}_{\tilde{w},\gamma}])\cdot {\rm Eul}(\qgot^{\gamma>0})= k[pt],\ k\geq 1$ in $H^*(\Fcal_\gamma,\Zbb)$.
\item $\sum_{\stackrel{\alpha\in\Rgot^+}{\langle\alpha,\gamma\rangle> 0}}\langle\alpha,\gamma\rangle=
\sum_{\stackrel{\tilde{\alpha}\in\tilde{\Rgot}^+}{\langle\tilde{\alpha},\tilde{w}_0\tilde{w}\gamma\rangle> 0}}\langle\tilde{\alpha},\tilde{w}_0\tilde{w}\gamma\rangle$.\end{itemize}
\end{theo}

\section{Example: the holomorphic Horn cone $\Horn(p,q)$}\label{sec:horn-p-q}

Let  $p\geq q\geq 1$. We consider the pseudo-unitary group $G=U(p,q)\subset GL_{p+q}(\Cbb)$ defined by the relation : $g\in U(p,q)$ if and only if 
$g{\rm Id}_{p,q}g^*={\rm Id}_{p,q}$ where ${\rm Id}_{p,q}$ is the diagonal matrice ${\rm Diag}({\rm Id}_p,-{\rm Id}_q)$.

We work with the maximal compact subgroup $K=U(p)\times U(q)\subset G$. We have the Cartan decomposition 
$\ggot=\kgot\oplus \pgot$, where $\pgot$ is identified with the vector space $M_{p,q}$ of $p\times q$ matrices through the map 
$$
X\in M_{p,q}\longmapsto\left(\begin{array}{cc} 0 & X \\ X^* & 0\\\end{array}\right).
$$
We work with the element $z_{p,q}=\frac{i}{2}{\rm Id}_{p,q}$ which belongs to the center of $\kgot$. The adjoint action of $z_{p,q}$ on  $\pgot$ corresponds to the standard complex structure on 
$M_{p,q}$.

The trace on $\mathfrak{gl}_{p+q}(\Cbb)$ defines an identification $\ggot\simeq \ggot^*=\hom(\ggot,\Rbb)$:  to $X\in\ggot$ we associate 
$\xi_X\in\ggot^*$ defined by $\langle \xi_X,Y\rangle=-{\rm Tr}(XY)$. Thus, the $G$-invariant cone $\Ccal_{G/K}$ defined by $z_{p,q}$ can be viewed as the following cone 
of $\ggot$:
$$
\Ccal(p,q)=\left\{X\in\ggot,\ {\rm Im}\left(\tr( gXg^{-1} {\rm Id}_{p,q})\right) \geq 0,\ \forall g\in U(p,q)\right\}.
$$

Let $T\subset U(p)\times U(q)$ be the maximal torus formed by the diagonal matrices. 
The Lie algebra $\tgot$ is identified with $\Rbb^p\times\Rbb^q$ through the map ${\bf d}:\Rbb^p\times\Rbb^q\to \ugot(p)\times\ugot(q)$ : 
${\bf d}_{x}={\rm Diag}(ix_1,\cdots,ix_p,i x_{p+1},\cdots, i x_{p+q})$. The Weyl chamber is
$$
\tgot_{\geq 0}=\left\{x\in \Rbb^{p}\times\Rbb^q,\ x_1\geq \cdots \geq x_p \ {\rm and}\  x_{p+1}\geq \cdots \geq x_{p+q}\right\}.
$$
Proposition \ref{prop:Gz} tells us that the $U(p,q)$ adjoint orbits in the interior of $\Ccal(p,q)$ are parametrized by the holomorphic chamber
$$
\Ccal_{p,q}=\left\{x\in \Rbb^p\times\Rbb^q, x_1\geq \cdots \geq x_p > x_{p+1}\geq \cdots \geq x_{p+q}\right\} \,\subset\, \tgot_{\geq 0}.
$$

\begin{defi}
The holomorphic Horn cone $\Horn(p,q):=\Horn^2(U(p,q))$ is defined by the relations  
$$
\Horn(p,q)=\left\{(A,B,C)\in(\Ccal_{p,q})^3,\ U(p,q){\bf d}_C\subset U(p,q){\bf d}_A+U(p,q){\bf d}_B\right\},
$$
\end{defi}

\medskip

Let us detail the description given of $\Horn(p,q)$ by Theorem B. For any $n\geq 1$, we consider the semigroup
$\wedge_n^+=\{(\lambda_1\geq \cdots\geq\lambda_n)\}\subset\Zbb^n$. If $\lambda=(\lambda',\lambda'')\in \wedge_p^+\times\wedge_q^+$, then
$V_\lambda:=V^{U(p)}_{\lambda'}\otimes V^{U(q)}_{\lambda''}$ denotes the irreducible representation of $U(p)\times U(q)$ with highest weight $\lambda$. 
We denote by ${\rm Sym}(M_{p,q})$ the symmetric algebra of $M_{p,q}$.

\begin{defi}\label{defi:horn-cone} 
\begin{enumerate}
\item $\horn^{\Zbb}(p,q)$ is the semigroup of $(\wedge_p^+\times\wedge_q^+)^3$ defined by the conditions:
$$
 (\lambda,\mu,\nu)\in \horn^{\Zbb}(p,q)\Longleftrightarrow \left[V_\nu\,:\,V_\lambda\otimes V_\mu\otimes {\rm Sym}(M_{p,q})\right]\neq 0.
$$
\item $\horn(p,q)$ is the convex cone of $(\tgot_{\geq 0})^3$ defined as the closure of 
$\Qbb^{>0}\cdot \horn^{\Zbb}(p,q)$.
\end{enumerate}
\end{defi}

Theorem B asserts that 
\begin{equation}
\Horn(p,q)=\horn(p,q)\,\bigcap \,(\Ccal_{p,q})^3.
\end{equation}

\medskip

In another article \cite{pep:horn-p-q}, we obtained a recursive description of the cones $\horn(p,q)$. 
This allows us to give the following description of the holomorphic Horn cone $\Horn(2,2)$.

\begin{exam}
An element $(A,B,C)\in(\Rbb^4)^3$ belongs to $\Horn(2,2)$ if and only if the following conditions holds:

\begin{equation*}
\boxed{
\begin{array}{rcl}
a_1 \geq a_2 & > & a_3 \geq a_4\\
b_1 \geq b_2 & > & b_3 \geq b_4\\
c_1 \geq c_2 & > & c_3 \geq c_4
\end{array}
}
\end{equation*}

$$
\boxed{a_1+a_2+a_3+a_4+b_1+b_2+b_3+b_4=c_1+c_2+c_3+c_4}
$$
$$
\boxed{a_1+a_2+b_1+b_2\leq c_1+c_2}
$$
\begin{equation*}
\boxed{
\begin{array}{rcl}
a_2+b_2&\leq & c_2\\
a_2+b_1&\leq & c_1\\
a_1+b_2&\leq & c_1
\end{array}
}
\end{equation*}
\begin{equation*}
\boxed{
\begin{array}{rcl}
a_3+b_3 &\geq & c_3\\
a_3+b_4 &\geq & c_4\\
a_4+b_3 &\geq & c_4
\end{array}
}
\end{equation*}
\begin{equation*}
\boxed{
\begin{array}{rcl}
a_2+a_4+b_2+b_4 &\leq & c_1+c_4 \\
a_2+a_4+b_2+b_4 &\leq & c_2+c_3\\
a_2+a_4+b_1+b_4 &\leq & c_1+c_3\\
a_1+a_4+b_2+b_4 &\leq & c_1+c_3\\
a_2+a_4+b_2+b_3 &\leq & c_1+c_3\\
a_2+a_3+b_2+b_4 &\leq & c_1+c_3
\end{array}
}
\end{equation*}
\end{exam}

\end{document}